\documentclass{amsart}
\usepackage{amsmath,amsfonts,amsthm,amssymb,nameref,mathtools,mathbbol,mathrsfs,amscd}
\usepackage[utf8]{inputenc}

\usepackage[margin=1.2in]{geometry}
\usepackage[utf8]{inputenc}

\usepackage{cite}
\usepackage{graphicx}
\usepackage[colorlinks=true, allcolors=blue]{hyperref}

\usepackage[all]{xy}
\usepackage{mathabx}

\usepackage{wasysym}
\usepackage{cancel,ulem}

\usepackage{tikz}

\usepackage{ifthen}

\newcommand{\showcomments}{yes}

\newsavebox{\commentbox}
{\ifthenelse{\equal{\showcomments}{yes}}%
{\footnotemark
    \begin{lrbox}{\commentbox}
    \begin{minipage}[t]{1.25in}\raggedright\sffamily\upshape\tiny
    \footnotemark[\arabic{footnote}]}%
{\begin{lrbox}{\commentbox}}}%
{\ifthenelse{\equal{\showcomments}{yes}}%
{\end{minipage}\end{lrbox}\marginpar{\usebox{\commentbox}}}%
{\end{lrbox}}}

\newcommand{\Isom}{\textrm{Isom}}

\newcommand{\revise}[1]{{#1}}

\newtheorem{theorem}{Theorem}[section]
\newtheorem{lemma}[theorem]{Lemma}
\newtheorem{proposition}[theorem]{Proposition}
\newtheorem{corollary}[theorem]{Corollary}

\newtheorem{question}[theorem]{Question}

\theoremstyle{definition}
\newtheorem{defn}[theorem]{Definition}

\theoremstyle{remark}
\newtheorem{remark}[theorem]{Remark}

\newtheorem*{claim}{Claim}

\newtheorem*{convention}{Convention}

\title[Uniform exponential growth for  proper product actions]{Uniform exponential growth for groups with proper product actions  on hyperbolic spaces}

\author{Renxing Wan}
\address{School of Mathematical sciences\\
East China Normal University\\
 Shanghai 200241, China P.R. }
\email{rxwan@math.ecnu.edu.cn}

\author{Wenyuan Yang}
\address{Beijing International Center for Mathematical Research\\
Peking University\\
 Beijing 100871, China P.R.}
\email{wyang@math.pku.edu.cn}

\keywords{Uniform exponential growth, product actions, hyperbolic spaces, product set growth, acylindrical actions}

\begin{document}

\begin{abstract}
This paper studies the locally uniform exponential growth and product set growth for a finitely generated  group $G$ acting properly on a finite product of hyperbolic spaces.  Under the assumption of coarsely dense orbits or shadowing property on factors, we prove that any  finitely generated non-virtually abelian  subgroup  has uniform exponential growth. These assumptions are   \revise{fulfilled} in many hierarchically hyperbolic groups, including mapping class groups, specially cubulated groups and BMW groups.

Moreover, if $G$ acts weakly acylindrically on each factor, we show that, with two exceptional classes of  subgroups, $G$ has uniform product set growth. As corollaries, this gives  a complete classification of subgroups with   product set growth for   any group acting discretely on a simply connected manifold with pinched negative curvature, for groups acting acylindrically   on trees, and for 3-manifold groups.
\end{abstract}
\maketitle

\section{Introduction}

Let $G$ be a finitely generated group and $S=S^{-1}$ be  a finite symmetric generating set. Denote by $\mathcal{G}(G, S)$ the Cayley graph of $G$ relative to $S$, equipped with word metric $d_S$. Let $S^{\leq n}$ denote the ball of radius   $ n$ in $\mathcal{G}(G, S)$ around the identity $1$. We often write  $|g|_S=d_S(1,g)$.

The \textit{exponential growth rate} of $G$ relative to $S$ is defined to be the following limit:
\begin{equation*}
\omega(G, S) = \lim_{n \rightarrow \infty}|S^{\leq n}|^{1/n}
\end{equation*}
which always exists by the
sub-multiplicative inequality   $|S^{\leq {n+m}}|\le |S^{\leq n}||S^{\leq m}|$.  If $\omega(G, S) > 1$, then   \revise{$\omega(G,S')>1 $ for any finite symmetric generating set $S'$} and we say that $G$ has \textit{exponential growth}.  Define the \textit{uniform exponential growth rate} of $G$ as
\begin{equation*}
\omega(G) = \inf \{\omega(G, S): \text{finite symmetric generating sets $S$} \}.
\end{equation*}
If $\omega(G) > 1$, we say that $G$ has \textit{uniform exponential growth}.


Many groups of exponential growth  have uniform exponential growth. This includes  non-elementary hyperbolic groups \cite{Kou}, relatively hyperbolic groups \cite{Xie}, solvable groups \cite{Alp}, \cite{Osi}, nontrivial amalgamated free products and HNN-extensions \cite{H}, one-relator groups \cite{BuH}, linear groups over a field of characteristic zero \cite{EMO}, and many hierarchically
hyperbolic groups \cite{ANSGP}. The first examples of groups   without
uniform exponential growth were constructed by Wilson \cite{W}. In  \cite{KS}, Kar-Sageev showed that if a group acts freely on a CAT(0) square complex, then either it has uniform exponential growth or it is virtually abelian. Later on, this is generalized to   any proper actions on CAT(0) square complexes by Gupta-Jankiewicz-Ng \cite{GJN}.

We say that $G$ has \textit{locally uniform exponential growth} if there is a constant $\omega_0 > 1$ such that every finitely generated non-virtually abelian subgroup $H$  has $\omega(H)\ge \omega_0$. This property is sometimes called ``uniform uniform exponential growth'' in literature. Such uniform bound on the growth of subgroups has been shown by Mangahas for the mapping class group \cite{Mah} (depending on the complexity of the surface), by Kar-Sageev for groups acting freely on square complexes \cite{KS}, and by Kropholler-Lyman-Ng for the automorphism group of a torsion-free one-ended hyperbolic group \cite{KLN}.

The goal of this paper is to study the uniform exponential growth for groups with a \textit{proper product action  on hyperbolic spaces}, i.e.: acting  properly on a finite product of hyperbolic spaces $X=\prod_{i=1}^lX_i$ endowed with $\ell_1$-metric and preserving each factor.  Any group acts properly on a proper hyperbolic space (e.g. the combinatorial horoball over the Cayley graph),  and the usual assumption to exclude this situation is requesting the existence of two independent loxodromic elements.  In this paper, we make use of   one of another two  assumptions. The first is that the product action has \textit{coarsely dense orbits}, which is equivalent to say that the product action is  co-bounded. The second is \textit{shadowing property} (see Definition \ref{ShadowPro}) on each factor. Roughly speaking, this ask for   a family of preferred quasi-geodesics in the Caylay graph that are projected to   un-parameterized quasi-geodesics in each factor $X_i$. This property is first studied by Masur-Minsky  in mapping class groups, where the preferred paths are given by hierarchy paths.

In  \cite{BBFa},  Bestvina-Bromberg-Fujiwara showed that the mapping class group $\mathrm{Mod}(\Sigma)$ of a finite type surface $\Sigma$ with $\chi(\Sigma)<0$ acts properly (with quasi-isometric embedded orbits) on a finite product of hyperbolic  graphs
\begin{align}\label{MCGAction}
\mathrm{Mod}(\Sigma) \curvearrowright \mathcal C(\mathbf Y^1)\times \cdots \times \mathcal C(\mathbf Y^l)
\end{align}
where $\mathcal C(\mathbf Y^i)$ denotes a quasi-tree of curve graphs of subsurfaces with the same ``color". The existence of finitely many colors (also called as BBF-coloring) for  all subsurfaces is due to \cite[Proposition 5.8]{BBFa}.  In \cite{BBFb},  $\mathcal C(\mathbf Y^i)$ can be further replaced with finite products of  quasi-trees (i.e. quasi-isometric to a simplicial tree). Proper action on finite product of quasi-trees is referred to as  \textit{Property (QT)} there. This property is characterized for 3-manifold groups in \cite{HTY}, and carefully examined in recent works of Button \cite{But19,But20}. The study of the so-called proper product actions  here on hyperbolic spaces  in this paper arise in a similar circle of ideas, for which we investigate on the uniform exponential growth.



\subsection{Uniform exponential growth}
Our first main result reads as follows.
\begin{theorem}[Theorem \ref{ueg}]\label{Thm1}
Suppose that a finitely generated group $G$ admits a proper product action  on hyperbolic spaces $X=\prod_{i=1}^lX_i$ with either coarsely dense orbits in $X$ or the shadowing property on each factor. Then $G$ has locally uniform exponential growth: any finitely generated non-virtually abelian subgroup has uniform exponential growth.
\end{theorem}

We remark that actually, if $G\curvearrowright X$ is co-bounded, then   the shadowing property holds    on each factor by Lemma \ref{cocompactshadowsfactor}.  This assumption is separated here for a better presentation.

Our first application is given to the class of hierarchically hyperbolic groups  (HHGs), introduced by Berhstock-Hagen-Sisto \cite{BHS19} with the prototype examples of mapping class groups and many cubical groups.  HHGs have been receiving very active research interests in recent years, with many examples found and techniques developed. See \cite{SisHHS} for a nice introduction to this topic.  In \cite[Prop. 3.10]{HP},   Hagen-Petyt showed that any HHG  with a \textit{BBF-coloring}   (in analog with the situation in $\mathrm{Mod}(\Sigma)$) admits a proper product action  on hyperbolic spaces with shadowing property on each factor. Note that  examples of HHGs which  do not virtually admit BBF-colorings were also produced in \cite{Hag21}. 

\begin{corollary}\label{HHG}
    Let $(G,\mathfrak G)$ be a  hierarchically hyperbolic group with a BBF-coloring. Then $G$ has locally uniform exponential growth.
\end{corollary}
It would be interesting to compare this  with a recent result of  Abbott-Ng-Spriano  \cite{ANSGP} that for a virtually torsion-free  HHG, either it has uniform exponential growth, or its Cayley graph is quasi-isometric to a direct product of the integers $\mathbb Z$ and a metric space. In the category of HHGs with BBF-colorings, our corollary   generalizes their result  with uniform exponential growth for any non-virtually abelian subgroups. On the other hand, the first examples of HHGs which are not virtually torsion-free were produced by Hughes in \cite{Hug}. Since the Hughes group is also a uniform lattice in product of trees, it has locally uniform exponential growth as a corollary of Theorem \ref{Thm1}.

Corollary \ref{HHG} also gives a partial answer to the question \cite[Question 1.6]{ANSGP} raised by Abbott-Ng-Spriano. Thus, their question is reduced to the following.
\begin{question}
    Does there exist a hierarchically hyperbolic group that is not colorable and does not have (locally) uniform exponential growth?
\end{question}

Another application is given to the class of Burger-Mozes-Wise groups. A group $G$ is a \textit{BMW group} if it acts by isometries on the product of two trees $T_1\times T_2$ such that every element preserves the product decomposition and the action on the vertex set of $T_1\times T_2$ is free and transitive.  Using these groups, Wise \cite{Wis} and Burger-Mozes \cite{BM} produced the first examples of non-residually finite and virtually simple CAT(0)
groups respectively. BMW groups have been extensively studied and have rich connections to the study of automata groups and commensurators (see Caprace's survey \cite{Cap}). A recent result of Abbott-Ng-Spriano \cite[Example 1.5]{ANSGP} gave a new proof of uniform exponential growth for  BMW groups. As a corollary of Theorem \ref{Thm1}, we have

\begin{corollary}
    BMW groups have locally uniform exponential growth.
\end{corollary}



\subsection{Product set growth} Recently, there is an increasing research interest in establishing {product set growth} in various classes of groups. This  is  actually a finer version of uniform exponential growth. We say that a subset $U$ of $G$ has \textit{product set growth} if there exist $\alpha, \beta>0$ such that
$|U^n|\geq (\alpha|U|)^{\beta n}$ for any $n\ge 1$ and $G$ has \textit{product set growth} if any finite generating subset $U$ (may not be symmetric) has product set growth for two uniform constants $\alpha,\beta>0$. If $U$ generates $G$ as a semi-group, then it implies $\omega(G,U)\ge (\alpha|U|)^\beta$, so the product set growth implies     uniform exponential growth. On the other hand, a group with infinite center clearly fails to have product set growth, but may have uniform exponential growth.

Such an inequality as above was first obtained for any non-elementary subgroup $\langle U\rangle$ in the class of hyperbolic groups in \cite{AL06}, where $U$ is assumed to be symmetric.     Until very recently, after  the work of non-symmetric case    in free groups \cite{Raz, Saf11}, establishing the product set growth for ``non-elementary" subgroups in  various classes of groups has becoming an active topic in Geometric Group Theory. Now, the product set growth has been known in  hyperbolic groups \cite{DS} and Burnside groups \cite{CS22}    for non-symmetric $U$, non-elementary subgroups of relatively hyperbolic groups \cite{But13, CJY}, groups acting acylindrically on trees \cite{CS21}, and mapping class groups \cite{Ker}  for symmetric $U$. In the recent work of Fujiwara-Sela \cite{FS23}, the product set growth serves as a first step to establish the well-ordered growth rates in hyperbolic groups (see \cite{Fuj23} for the results on more general classes of groups).

\begin{convention}
In practice, it appears  harder to prove the  product set growth for a non-symmetric subset $U$. In this paper, we always consider the case  that (possibly non-symmetric) $U$ generates a group denoted by $\langle U\rangle$ as \textbf{semigroup}:   every element in $\langle U\rangle$ can be written as a word over the alphabet set $U$ (and of course $U\cup U^{-1}$). To be brief, we shall say that $U$ generates $\langle U\rangle$ as a semigroup.
\end{convention}


The next goal   of this paper is aiming to  explore the   product set growth for groups with proper product action  on hyperbolic spaces $X=\prod_{i=1}^lX_i$ with  \textit{weakly acylindrical} action on each factor $X_i$.  Generalizing the usual acylindrical one,  a  weakly acylindrical action is defined by Delzant \cite{TD} as a uniformized version of WWPD action in \cite{BBFa}. Roughly speaking, the projection of the quasi-axis for any loxodromic element $g\in G$ to any other $G$-translate has a linear upper bound comparable to the stable length of $g$. See Definition \ref{WACYLDef} for details.

Before stating the general form of the main result, Theorem \ref{Thm2}, we would like to first point out some corollaries of it, where there is only one factor, i.e. $l=1$.

Thanks to Margulis Lemma, Delzant observed in \cite[Example 1]{TD} that any discrete group action on a  simply connected, complete, and negatively pinched Riemannian manifold is weakly acylindrical.  As  a by-product to the proof of Theorem \ref{Thm2}, we are able to prove the product set growth for such groups.  This strengthens the earlier result obtained by  Besson-Courtois-Gallot  \cite{BCG} \revise{that non-elementary groups acting properly on a Cartan-Hadamard manifold with negatively pinched curvature have uniform exponential growth}.

\begin{theorem}[Theorem \ref{RiemPSGThm2}]\label{RiemPSGThm}
    Let $M$ be a Riemannian manifold with pinched negative curvature. Then there exist constants $\alpha,\beta>0$ depending only on dimension and curvature such that for any finite (possibly non-symmetric) $U\subset \pi_1(M)$ generating a group $\langle U\rangle$ as  a semi-group,   one of the following is true:
    \begin{enumerate}
    \item $\langle U\rangle$ is virtually nilpotent.
    \item $|U^n|\geq (\alpha|U|)^{\beta n}$ for every $n\in \mathbb N$.
\end{enumerate}
\end{theorem}

Recall that a subgroup $H$ of a relatively hyperbolic group is called elementary if the limit set of $H$ contains no more than two points in Bowditch boundary. The proof of Theorem \ref{Thm2}, together with Lemma \ref{ShortHypLem}, also improves the result in \cite{CJY}.
\begin{theorem}\label{RHGPSGThm}
Let $G$ be a non-elementary relatively hyperbolic group. Then there exist constants $\alpha,\beta>0$ with the following property. Let $U$ be any finite (possibly non-symmetric) subset such that $U$ generates a non-elementary subgroup $\langle U\rangle$ as a semi-group. Then  $|U^n|\geq (\alpha|U|)^{\beta n}$ for every $n\in \mathbb N$.
\end{theorem}

The following result removes the symmetry of $U$ in a result of  Cerocchi-Sambusetti  \cite[Theorem 1]{CS21}. \revise{This could be also deduced  from  Delzant-Steenbock \cite[Theorem 1.11]{DS} \footnote{We thank Alice Kerr to point out this to us.}, where the optimal $\beta =1/2$ is actually obtained.   We record the result for further reference.} 
\begin{theorem}\label{AcylTreeThm}
Let $G$ be any group acting acylindrically  on a simplicial tree $X$. Then there exist constants $\alpha,\beta>0$ with the following property. Let $U$ be any finite possibly non-symmetric subset such that $U$ generates as a semi-group a non-virtually cyclic subgroup $\langle U\rangle$ without global fixed point in $X$. Then  $|U^n|\geq (\alpha|U|)^{\beta n}$ for every $n\in \mathbb N$.
\end{theorem}

From the above two results, we can classify subgroups with the product set growth in 3-manifold groups. \revise{A compact 3-manifold $M$  is called a mixed manifold (resp. graph manifold) if  $M$ contains (resp. no) hyperbolic pieces in its nontrivial JSJ decomposition. Accordingly, the fundamental group is  relatively hyperbolic  or acts   acylindrically  on a simplicial tree.} If a compact 3-manifold $M$ has trivial JSJ decomposition, then the Geometrization Theorem tells us that $M$ supports one of the eight Thurston geometries.  Except the hyperbolic and Sol geometry, the other six geometries admit a central extension of surface subgroups, so it is straightforward to describe subgroups with product set growth.   The question whether the product set growth holds in (non-nilpotent) solvable groups seems widely open. With Sol geometry undetermined yet, we can derive the product set growth for  all the remaining cases from Theorems \ref{RHGPSGThm} and \ref{AcylTreeThm}.
\begin{corollary}\label{3MfdPSGCor}
Let $M$ be a connected, compact, oriented and irreducible 3-manifold with empty or tori boundary, either with hyperbolic geometry or with nontrivial JSJ decomposition and that does not support the Sol geometry. Then there exist constants $\alpha,\beta>0$ such that for any  (possibly non-symmetric) finite subset  $U\subset \pi_1(M)$, one of the following is true:
\begin{enumerate}
\item $\langle U\rangle$ is virtually abelian of rank at most 2 or a subgroup of a Seifert sub-manifold group with infinite center. 
\item $|U^n|\geq (\alpha|U|)^{\beta n}$ for every $n\in \mathbb N$.
\end{enumerate}
\end{corollary}

Now let us state our main result about product set growth for a proper product action with weakly acylindrical actions on factors. For a loxodromic element $g$ on a hyperbolic space, we denote $E(g)$ as  the  subgroup of elements which fix setwise the endpoints of the quasi-axes $L_g$ of $g$.
\begin{theorem}[Theorem \ref{PSG}]\label{Thm2}
Let $G$ admit  a proper product  action on    hyperbolic spaces $X=\prod_{i=1}^lX_i$,  with shadowing property on factors. Suppose that $G$ acts weakly acylindrically on each factor. Then there exist $\alpha,\beta >0$ depending on $\delta,l$ such that for every finite symmetric subset $U\subset G$, at least one of the following  is true:
\begin{enumerate}
    \item $\langle U\rangle$ is virtually abelian.
    \item There exists a factor $X_i$ so that no uniform bound exists on $|E(h)\cap E(b)\cap \langle U\rangle|$ for any two independent loxodromic elements $b,h\in \langle U\rangle$ on $X_i$.
    \item $|U^n|\geq (\alpha|U|)^{\beta n}$ for every $n\in \mathbb N$.
\end{enumerate}
Moreover, if there is no lineal action of $\langle U\rangle$ on factors, then $U$ can be assumed non-symmetric in the above statement.
\end{theorem}

At last, we would like mention an    application of Theorem \ref{Thm2} concerning   mapping class groups $\mathrm{Mod}(\Sigma)$ of a surface $\Sigma$. See Corollary \ref{MCGPSG} for the more elaborated version where $\langle U\rangle$ is a subgroup. 

\begin{corollary}\label{MCGPSGCor}
    Let $\Sigma$ be an orientable finite-type surface of $\chi(\Sigma)<0$ with possibly cusps   without boundary. Then there exist $\alpha, \beta>0$ depending only on $\Sigma$ such that for every finite possibly non-symmetric subset $U$ generating $\mathrm{Mod}(\Sigma)$ as a semi-group,  one of the following is true:
\begin{enumerate}
    \item There exist a factor $\mathcal C(\mathbf Y^i)$ in (\ref{MCGAction}) and two independent loxodromic elements $b,h$ on $\mathcal C(\mathbf Y^i)$ such that $|E(h)\cap E(b)|$ is infinite.
    \item $|U^n|\geq (\alpha|U|)^{\beta n}$ for every $n\in \mathbb N$.
\end{enumerate}
\end{corollary}

In a recent work of Kerr,  the uniform product set growth has been classified for mapping class groups  \cite[Theorem 1.0.2]{Ker}, where $U$ is assumed to be symmetric: $\langle U\rangle$  has product set growth if and only if $Z(H)$ is finite for any finite-index subgroup $H\le \langle U\rangle$.  The question whether    the symmetry of $U$ could be removed remains open.  On one hand, Corollary \ref{MCGPSGCor} allows $U$ to be non-symmetric; on the other hand, it is not clear whether the item (1) happens exactly when $\langle U\rangle$ have virtually infinite center.   


\subsubsection*{\textbf{Structure of the paper.}} This paper is organized as follows. In Section \ref{PreSection}, we recall some standard materials in Gromov's hyperbolic geometry and give definitions of shadowing property and weakly acylindrical actions. In addition, we removes the symmetry condition in a lemma of Cui-Jiang-Yang \cite{CJY} which produces a short loxodromic element from a large displacement. In Section \ref{ShortLoxoSection}, we utilise the  shadowing property on each factor to transfer a large displacement from the product space to some factor space. In Section \ref{UEGSection}, we first show that any non-virtually abelian group has an arbitrarily large commutator subset, which  in turn is used to prove the existence of two short independent loxodromic elements on some factor space. The proof of Theorem \ref{Thm1} is then finished with a result of Breuillard-Fujiwara \cite{BF18}. In Section \ref{PSGSection}, we first give some preliminary results about product set growth  for non-symmetric subsets. In \textsection \ref{SSubPSGOnOneFactor}, we present detailedly the method used in \cite{CJY} to prove the product set growth on one factor, Theorem \ref{PSGonOneFactor}. This  is an essential step to get the main result Theorem \ref{Thm2} on a finite product of hyperbolic spaces, whose proof is then completed in \textsection\ref{SSubProofPSG}. As corollaries, we give some applications including Theorems \ref{RHGPSGThm}, \ref{AcylTreeThm} and Corollaries \ref{3MfdPSGCor}, \ref{MCGPSGCor} about relatively hyperbolic groups, groups acting acylindrically on trees, 3-manifold groups and mapping class groups  in \textsection\ref{SSubApps}. At last, in Section \ref{SSubHadamard},  the proof of Theorem \ref{RiemPSGThm} is presented for the uniform product set growth for the fundamental group of a Riemannian manifold with pinched negative curvature.

\subsection*{Acknowledgments}
We are grateful to Prof. Thomas Delzant for helpful communications, in particular, for providing the proof  of Lemma \ref{WeaklyAcy}. {Thanks go to Sam Hughes and  Alice Kerr for pointing out the references  and  for several  very helpful comments. We greatly thank the referee for many helpful comments and criticism  on the writing, which helped much improve the exposition.}    W. Y. is supported by National Key R \& D Program of China (SQ2020YFA070059), and   National Natural Science Foundation of China (No.  12131009 and No. 12326601).

\section{Preliminary}\label{PreSection}

In this section, we   introduce the preliminaries that we will use later on.

\subsection{Gromov hyperbolic spaces and classification of isometric actions}\label{GromovHypSpace}

Let $(X,|\cdot|)$ be a geodesic metric space. For any $a,b \in X$, we often denote by $[a, b]$ a choice of a geodesic between $a, b$, if there is no ambiguity in context. For any $A,B\subset X$, we denote $d(A,B):=\inf_{a\in A, b\in B}|a-b|$. A geodesic triangle $\Delta=\Delta(a,b,c)$ consists of three geodesics $[a,b],[b,c],[c,a]$.  For any $a,b,c\in X$, define  \textit{Gromov product} $\langle a,c\rangle_b$    as follows:
\begin{equation*}
\langle a,c\rangle_b=\frac{|a-b|+|c-b|-|a-c|}{2}.
\end{equation*}

We say that $X$ is a \textit{$\delta$-hyperbolic space} for some $\delta\geq 0$ if every geodesic triangle has a point  within a distance $\delta$ to each   side.




By enlarging $\delta$ by a constant times,  one can show the following stronger thin-triangle property:
\begin{itemize}
    \item If $x\in [b,a]$ and $y\in [b,c]$ so that $|b-x|=|b-y|\le \langle a,c\rangle_b$, then $|x-y|\le \delta$.
\end{itemize}

Unless otherwise stated, all hyperbolic spaces in this paper are assumed to be $\delta$-hyperbolic  for a fixed constant $\delta$.

\begin{defn}
Given  $k\geq 1, c>0$, a map between two metric spaces $f: (X,|\cdot|_X)\rightarrow (Y,|\cdot|_Y)$ is called a \textit{$(k,c)$-quasi-isometric embedding} if the following inequality holds
\begin{equation*}
\frac{|x-x'|_X}{k}-c\leq |f(x)-f(x')|_Y\leq k|x-x'|_X+c,
\end{equation*}
for all $x,x'\in X$. Furthermore, if there exists $R>0$ such that $Y\subset N_R(f(X))$, then $f$ is called a \textit{$(k,c)$-quasi-isometry}.
\end{defn}

When the values of $k,c$ are clear in context or do not matter, we omit them and just say $f$ is a quasi-isometric embedding or quasi-isometry.

If a path $p$ in $X$ is  rectifiable, denote by $Len(p)$ the length of $p$, and by $p_-,p_+$ the endpoints of $p$.
\begin{defn}
A path $p$ in $(X,d)$ is called a \textit{$(k,c)$-quasi-geodesic} for $k\geq 1,c\geq 0$ if the following holds
\begin{equation*}
Len(q)\leq k|q_--q_+|+c
\end{equation*}
for any rectifiable subpath $q$ of $p$.
\end{defn}
Sometimes, a $(k,k)$-quasi-geodesic will be called a $k$-quasi-geodesic for simplicity.

\begin{lemma}[Morse Lemma]\label{morse}
Let $p,q$ be two $(k,c)$-quasi-geodesics with the same endpoints in a hyperbolic space $X$. Then there exists $D=D(\delta,k,c)$ such that  $p\subset N_D(q)$ and $q\subset N_D(p)$.
\end{lemma}

We shall need the following (un)parametrized version of quasi-geodesic paths.

\begin{defn}
Let $I$ be a possibly infinite interval of $\mathbb R$.
A map $p$ from   $I$   to a metric space $X$ is called an \textit{unparametrized $(k,c)$-quasi-geodesic} for $k\ge 1, c\ge 0$ if there exists a monotone non-decreasing map $\tau: I\to I$ such that $p\cdot \tau$ is a $(k,c)$-quasi-geodesic in the usual sense.
\end{defn}

Another useful fact in hyperbolic geometry is the following lemma.


\begin{lemma}\label{projection}
Let $\alpha$ be a $(k,c)$-quasi-geodesic in a hyperbolic space $X$ and $q: X\to \alpha$ be the closest point projection. Then there exists a constant $C=C(k,c, \delta)$ so that  for any $x,y\in X$,
\begin{equation*}
    |q(x)-q(y)|\leq |x-y|+C
\end{equation*}
\end{lemma}
\revise{
\begin{defn}\cite[Definition 3.12]{BH}
Let $X$ be a  $\delta$-hyperbolic space with a basepoint $o$. A sequence $(x_n)_{n=1}^{\infty}$ in $X$ \textit{converges at infinity} if $\langle x_i,x_j\rangle_o\to \infty$ as $i,j\to \infty$.
The Gromov boundary $\partial X$ is defined as the following set, equipped with an appropriate metrizable topology,
$$\partial X:=\{(x_n)\mid(x_n)\text{ converges at infinity}\}/\sim$$
where $(x_n)\sim (y_n)$ if $\langle x_i,y_j\rangle_o\to\infty$ as $i,j\to \infty$.
\end{defn}
}

By Gromov \cite{G}, the isometries of a hyperbolic space $X$ can be subdivided into three classes. A nontrivial element $g\in \Isom(X)$ is called \textit{elliptic} if some $\langle g\rangle$-orbit is bounded. Otherwise, it is called \textit{loxodromic} (resp. \textit{parabolic}) if it has exactly two fixed points (resp. one fixed point)  in the Gromov boundary $\partial X$ of $X$. Note that, if an isometry has three fixed points in the boundary, then it leaves the center coarsely invariant, so it must be an elliptic element. Hence, this gives a complete classification of isometries.

If $g$ is a loxodromic element, any $\langle g\rangle$-invariant quasi-geodesic between the two fixed points shall be refereed as    a \textit{quasi-axis} for $g$. \revise{Two loxodromic elements are called \textit{independent} if their fixed points in $\partial X$ are disjoint.} Loxodromic elements play an important role in what follows. We list a few well-known facts for later use.

\begin{lemma}\cite[Lemma 9.2.2]{CDP}\label{LoxoCriterionLem}
If an isometry $g$ on a hyperbolic space $X$ satisfies  $$|o-go| \ge 2\langle o, g^2 o\rangle_{go} + 6\delta$$ for some point $o \in X$, then $g$ is loxodromic.
\end{lemma}

Let  $g$ be an isometry on a metric space $(X,|\cdot|)$. The \textit{stable translation length} of $g$ is defined to be $$\tau(g)=\lim_{n\to \infty}\frac{|x-g^nx|}{n}$$ for any $x\in X$, and the \textit{minimal translation length} of $g$ is defined to be $$\lambda(g)=\inf_{x\in X}|x-gx|.$$

\begin{lemma}\cite[Prop. 10. 6.4]{CDP}\label{StableTransLength}
    If $g$ is a loxodromic element on a hyperbolic space $X$, then $|\tau(g)-\lambda(g)|\le 16\delta$.
\end{lemma}

\begin{lemma}\cite[Lemma 2.4]{CJY}\label{LoxoCriterion2Lem}
If $g, h$ are two  isometries satisfying  $$\frac{1}{4}\min\{|go-o|, |ho- o|\} \ge L\ge \max\{\langle go, h^{-1} o\rangle_o,\langle g^{-1}o, h o\rangle_o\} + \delta$$ for some point $o \in X$ and some $L>0$. Then
\begin{enumerate}
    \item
    $gh$ is loxodromic.
    \item
    there  \revise{exists a constant $ c>0$} depending only on $L$ and $\delta$ such that the concatenated path $\alpha:=\bigcup_{i\in\mathbb Z} (gh)^{i}\left([o, go]\cdot g[o,ho]\right)$ is a   \revise{$c$}-quasi-geodesic.
\end{enumerate}
\end{lemma}

Furthermore, an unbounded isometric group action on Gromov hyperbolic spaces can be classified into the following four types:
\begin{enumerate}
    \item \textit{horocyclic} if it has no loxodromic element;
    \item \textit{lineal} if it has a loxodromic element and any two loxodromic elements have the same fixed points in Gromov boundary;
    \item \textit{focal} if it has a loxodromic element, is not lineal and any two loxodromic elements have one common fixed point;
    \item \textit{of general type} if it has two loxodromic elements with no common fixed point.
\end{enumerate}

Let $G\curvearrowright X$ be an unbounded isometric action on a hyperbolic space. A point $p$ in $\partial X$ is called a \textit{limit point} of $G$ if there exists a sequence of elements $g_n\in G$ such that $g_n o\to p$ for some (or any) $o\in X$.
The set of all limit points forms the \textit{limit set} denoted by $\Lambda G$. \revise{By definition, the limit set of a lineal action consists of two distinct points. By restricting to the union of geodesics between the two limit points, we may assume that a lineal action $G$ acts on a quasi-line. } Recall that  a \textit{quasi-tree} (resp. \textit{quasi-line}) is a graph quasi-isometric to a simplicial tree (resp. real line).




A lineal action of a group $G$ on a hyperbolic space is intimately related to the following notion.

\begin{defn}\cite[Definition 2.4]{M}
A map $f: G\rightarrow \mathbb R$ is called a \textit{quasi-character} if there exists $K\ge 0$ such that for any $g,h\in G$ , $$|\Delta f(g,h)|=|f(gh)-f(g)-f(h)|\leq K$$

If in addition, $f(g^n)=nf(g)$ for all $n\in \mathbb Z$ and $g\in G$, then $f$ is a \textit{pseudo-character}. The quantity $\Delta(f):=\sup\{|\Delta f(g,h)|:g,h\in G \}$ is called \textit{defect} of $f$.
\end{defn}

\begin{remark}\label{pcha}
It is well-known that the homogenuous version of a quasi-character defined as follows
$$\varphi(g)=\lim_{n\rightarrow \infty}\frac{f(g^n)}{n}$$
is a pseudo-character, so that $\|f-\varphi\|_\infty$ bounded above by $\Delta(f)$.
A pseudo-character $\varphi$ takes the constant on conjugacy classes: $\varphi(g)=\varphi(aga^{-1})$ for any $a, g\in G$. This implies that $|\varphi([a,b])|\le \Delta \varphi\le 4\Delta f$ for any $a, b\in G$.
\end{remark}

A lineal action of a group $G$ on a hyperbolic space $X$ is called \textit{orientable} if no element of $G$ permutes the two limit points of $G$ on $\partial X$. \revise{It is clear that an index at most 2 subgroup of $G$  is  orientable.}

\begin{lemma}
Let $G\curvearrowright X$ be an orientable  lineal  action on a quasi-line. Then $G$ admits a pseudo-character $\varphi: G\to \mathbb R$ so that the induced quasi-action is quasi-conjugated to the action $G\curvearrowright X$.
\end{lemma}
\begin{proof}
The proof is well-known; we present it for completeness.

By \cite[Lemma 3.7]{M}, there exists  an almost isometry $p: X\to \mathbb R$ with almost isometry $q$ as quasi-inverse. Define  a quasi-action of $G$ on $\mathbb R$ by $g\circ y =p(g\cdot qy)$ for any $y\in \mathbb R$. Thus, each $g\in G$ acts on  $\mathbb R$ by almost translation: there exists a uniform $D>0$ independent of $g$ and $y$ so that $|g\circ y-\tau_g(y)|\le D$ where $\tau_g$ is a translation on $\mathbb R$. Moreover, $p$ quasi-conjugates the action $G\curvearrowright X$ to the quasi-action  of $G$ on $\mathbb R$, i.e.: there exists $C>0$ so that $|p(g\cdot x)-g\circ p(x)|\le C$ for all $x\in X$.  Fix any $x\in X$.  Define $\varphi: G\to\mathbb R$ by $\varphi(g):=p(g\cdot x)$. Then $|\varphi(gh)-gh\circ (px)|\le |p(gh\cdot x) - gh\circ (px)|\le C$. As $g, h$ act  by almost translations on $\mathbb R$, we have $|\varphi(gh)-\varphi(g)-\varphi(h)|\le 3C+3D$.
\end{proof}

As an immediate consequence, we have the following result.

\begin{lemma}\label{bounded}
Let $G\curvearrowright X$ be an orientable lineal action on a quasi-line. Then there exists $C=C(\delta)>0$ such that for any $g,h\in G$ and $o\in X$, $|[g,h]o-o|\leq C$.
\end{lemma}

\subsection{Shadowing property}

Now we define \textit{shadowing property} for an isometric action on a metric space.




\begin{defn}\label{ShadowPro}
We say that an isometric action of $G$ on a  metric space $X$ has \textit{shadowing property} if for some (or any) finite symmetric generating set $F$ of $G$ and a basepoint $o\in X$, the Cayley graph $\mathcal{G}(G,F)$  admits a path system $\mathcal{P}$  of  uniform quasi-geodesics with the following properties:
\begin{enumerate}
    \item Any two elements $g,h\in G$ are connected by some  $\gamma \in \mathcal{P}$.
    \item The (connected) subpaths of any $\gamma\in \mathcal{P}$ are contained in $\mathcal{P}$.
    \item For any $\gamma\in \mathcal{P}$, the image path $\tilde{\gamma}=\Phi(\gamma)$ in $X$ is a uniform un-parameterized quasi-geodesic, where  $\Phi: \mathcal G(G,F)\to X$ is the orbit map given by $g\mapsto g\cdot o$.

\end{enumerate}
\end{defn}
\begin{remark}
The shadowing property does not depend on the choice of finite generating set $F$ and the basepoint $o\in X$. Indeed, the Cayley graphs for different generating sets are bi-Lipschitz, so the path system $\mathcal P$ with properties (1) (2) for one generating set is a  path system of quasi-geodesics with different parameters and with same properties (1) (2). Moreover, if an orbital image path $\tilde{\gamma}$  is a quasi-geodesic for one basepoint $o$, then it is true for any other basepoint $o'\in X$  (with different parameters).
\end{remark}

By Lemma \ref{morse}, an immediate consequence of the above property (3) is as follows.
\begin{lemma}\label{shadow}
Assume that $G\curvearrowright X$ has shadowing property, where $X$ is hyperbolic. There exists a constant $D>0$ depending only on the path system $\mathcal P$ with the following property.
For any $g,h\in G$, let $\gamma$ be a path connecting $g$ to $h$ in $\mathcal P$, then
\begin{equation*}
[go, ho]\subset N_D(\Phi(\gamma)),\;\; \Phi(\gamma)\subset N_D([go,ho]).
\end{equation*}
\end{lemma}


The shadowing property applied to hyperbolic spaces implies the cobounded action.
\begin{lemma}\label{ShadowCoboundedLem}
Assume that  $G\curvearrowright
 X$ has shadowing property, where $X$ is hyperbolic. Then $G$ acts on $X$ with quasi-convex $G$-orbit. In particular, by taking the weakly convex hull of any $G$-orbit, we can assume that $G$ acts co-boundedly on $X$.
\end{lemma}
\begin{proof}
Let $o\in X$ be a basepoint. We shall show that $Go$ is a quasiconvex subset of $X$. It suffices to find a constant $R>0$ such that any geodesic with endpoint in $Go$ is contained in $N_R(Go)$. Consider a geodesic  $\alpha=[go, ho]$ in $X$, and a quasi-geodesic $\gamma\in \mathcal P$ connecting $g$ to $h$ given by the definition of shadowing property. It is clear that the image of $\gamma$ under the orbital map is contained in $Go$. As $\Phi(\gamma)$ is a unparametrized quasi-geodesic, the Morse Lemma shows that  $\alpha$ is contained in $N_R(Go)$.

Let $\tilde X$ be the weakly convex hull of $Go$, which is the union of all geodesics with two endpoints in $Go$. By hyperbolicity, $\tilde X$ is contained in a uniform neighbourhood of $Go$, on which $G$ acts co-boundedly. It is starightforward  that $\tilde X$  with induced length metric is quasi-isometrically embedded into $X$. The conclusion follows.
\end{proof}

\revise{Let  $G\curvearrowright X=\prod_{i=1}^lX_i$ be a product action. We say that $G\leq \text{Isom}(X)$ \textit{preserves each factor} if the embedding $G\to \text{Isom}(X)$ factors through $\prod_{i=1}^l\text{Isom}(X_i)$. In this case, $G$ acts  isometrically   on each $X_i$,  with possibly a large kernel. }

Here is a useful fact to obtain shadowing property in the following setup.
\begin{lemma}\label{cocompactshadowsfactor}
Suppose that a finitely generated group $G$ acts properly and \revise{co-boundedly} on a finite product $X=\prod_{i=1}^l X_i$ of geodesic metric spaces and preserves the factors. Then $G\curvearrowright X_i$ has the shadowing property for each $1\le i\le l$.
\end{lemma}
\begin{proof}
\revise{Assume that $X$ is equipped with either  $\ell_1$-metric or $\ell_2$-metric. By Milnor-Svarc Lemma, $G$ admits a finite symmetric generating set $F$ so that the Cayley graph $\mathcal{G}(G,F)$ is quasi-isometric to $X$. To be precise, let $\phi:\mathcal{G}(G,F)\to X$ be the $G$-equivariant quasi-isometry and $\psi: X\to \mathcal{G}(G,F)$ be its $G$-equivariant quasi-inverse.

In order to achieve the shadowing property in Definition \ref{ShadowPro}, we first realize  a similar   path system $\mathcal Q$ on $X$ so that any two points $x,y$ are connected by one and only one $\ell_1$ (resp. $\ell_2$)-geodesic in $\mathcal Q$. Moreover, the path system $\mathcal Q$ is closed under taking subpaths, and is invariant under the $G$--action.

We now explain the construction of $\mathcal Q$. It is based on the following property for the $i$-th projection of $X\to X_i$ for each factor $X_i$ :
\begin{itemize}
    \item Any $\ell_1$ (resp. $\ell_2$)-geodesic $\gamma$ projects to an unparametrized geodesic $\gamma_i$ in $X_i$. The same holds true for any subsegment of $\gamma$.
\end{itemize}
The $\ell_2$-metric on $X$ is uniquely geodesic, so the construction of $\mathcal Q$ is obvious. For $\ell_1$-metric, define $d(x,y)=\sum_{i=1}^l d_i(x_i,y_i)$ for $(x_i), (y_i)\in X$. The geodesic between $x$ and $y$ is a concatenation of  segments $$\bigcup_{1\le i\le l} [(y_1,\ldots,y_{i-1},x_i,x_{i+1},\ldots, x_l), (y_1,\ldots, y_i,x_{i+1},\ldots,x_l)]$$ where the $i$-th segment only varies on the coordinate in $X_i$.  We  can order the factors $X_i$ in $X$ with the following property. Let $\mathcal Q$ be the set of  $\ell_1$-geodesics $\gamma$  among the geodesics from $x$ to $y$ so that the $l$-tuple $(d_i(x_i,y_i))_{1\le i\le l}$ is lexicographical minimal. Any subpath of $\gamma$ is lexicographical minimal, so $\mathcal Q$ is closed by taking subpaths.  Moreover, $\mathcal Q$ is  invariant under the factor-preserving isometric action of $G$.

Therefore,    setting $\mathcal P:=\psi(\mathcal Q)$ gives the desired path system on $\mathcal{G}(G,F)$, where $\phi: \mathcal{G}(G,F)\to X$ is the above  $G$-equivariant  quasi-isometry. }
\end{proof}
\subsection{Weakly acylindrical action}
Assume that $G$ acts isometrically on a hyperbolic space $X$.

For any $x\in X$, the \textit{displacement} of a subset $S\subset G$ at $x$ is defined as follows:
\revise{
\begin{equation*}
    \lambda(S, x):=\max_{s\in S}|x-sx|.
\end{equation*}}

Define the \textit{displacement} of $S$ on a subset $A\subset X$ as
\begin{equation}\label{displacement}
    \lambda(S, A):=\inf_{x\in A} \lambda(S, x)=\inf_{x\in A}\max_{s\in S}|x-sx|.
\end{equation}
If $S=\{g\}$ consists of one element, we also write $\lambda(g, x)$ and $\lambda(g, A)$.

\begin{defn}
The action of $G$ on a hyperbolic space $X$ is \textit{acylindrical} if there exist an integer $N$ and a positive number $K$ such that for every pair of points $x,y$ with $|x-y|\ge K$, the set $\{g\in G\mid |x-gx|+|y-gy|\le 1000\delta\}$ is finite with cardinality $\le N$.
\end{defn}

\revise{
If $g$ is a loxodromic element with $\lambda(g,X)>1000\delta$, then \cite[Lemma 11.3]{BF18} gives a $\langle g\rangle$-invariant $(\frac{9}{10},24\delta)$-quasi-geodesic $L_g$ in $X$. Let $E(g)$  be  the  subgroup of elements $h\in G$ which preserves  the two endpoints of the quasi-axes $L_g$.}

\begin{defn}\cite[Definition 6]{TD}\label{WACYLDef}
The action of $G$ on a hyperbolic space $X$ is \textit{weakly acylindrical} if there exists a constant $D>0$ such that for every loxodromic element $g$, and every $h\notin E(g)$, the diameter of the projection of $hL_g$ to $L_g$ is bounded by $D(1+\lambda(g, X))$.

If the constant $D$ is specified, the action is called $D$-\textit{weakly acylindrical}.
\end{defn}

\begin{remark}\label{AcyImplyWeaklyAcy}
    \revise{In \cite[Proposition 5]{TD}, Delzant showed that an acylindrical action on a hyperbolic space is weakly acylindrical. He further showed that any proper action on Hadamard manifold with negatively pinched curvature is weakly acylindrical, as recalled in Theorem \ref{WeaklyAcy}. This is an important fact  we make use of to prove Theorem \ref{RiemPSGThm}.}
\end{remark}

\begin{lemma}\label{alternative}
Suppose that  $G$ acts weakly acylindrically on a hyperbolic space $X$ with a loxodromic element. Then $G\curvearrowright X$ is either lineal or of general type.
\end{lemma}

\begin{proof}
Suppose that $G\curvearrowright X$ is focal. Then   any two loxodromic elements $g,h\in G$ have a common fixed point on the Gromov boundary $\partial X$, so the projection of $hL_g$ to $L_g$ is unbounded. This contradicts  the definition of weakly acylindrical actions. The conclusion thus follows by the classification of isometric actions on hyperbolic spaces.
\end{proof}


\begin{lemma}\label{LinearBILem}
Assume that $G$ acts weakly acylindrically on a hyperbolic space $X$. Let $h$ be a loxodromic element in $G$ so that for some point $o\in X$ and $c>0$,  the path  $\gamma=\bigcup_{n\in \mathbb Z} [h^no, h^{n+1}o]$ is a $c$-quasi-geodesic  in $X$.  Then for any given $\theta>0$,  there exists $D=D(c,\delta,  \theta)$ independent of the point $o$ such that for any $f\notin E(h)$,  we have $$\text{diam}(\gamma\cap N_R(f\gamma))\le D\cdot(  |o- ho|+1)$$
where $R:=\theta\cdot |o- ho|.$
\end{lemma}
\begin{proof}
Let $\pi:X\to \gamma$ be a closest projection map, assigning a point $x\in X$ to a shortest projection point $\pi(x)$ on $\gamma$. It follows from the definition of weakly acylindrical action and Morse Lemma that there exists a constant $D_0=D_0(c,\delta)$ such that $\text{diam}(\pi(f\gamma))\leq D_0(1+\lambda(h))$. Set $p:=\pi(f\gamma)\subset \gamma$. Then $\text{diam}(p)\leq D_0(1+\lambda(h))\leq D_0(1+|o-ho|)$.

For any $x\in \gamma\cap N_R(f\gamma)$, there exists $y\in f\gamma$ with $|x-y|\leq R$. It follows that $|y-\pi(y)| \leq |y-x|\leq R$ and then $|x-\pi(y)|\leq |x-y|+|y-\pi(y)|\leq 2R$. As a result, $d(x,p)\leq 2R$ and then $\gamma\cap N_R(f\gamma)\subset N_{2R}(p)$.

Thus, $\text{diam}(\gamma\cap N_R(f\gamma))\leq \text{diam}(p)+4R\leq (D_0+4\theta)|o-ho|+D_0$. By setting $D=D_0+4\theta$, we complete the proof.
\end{proof}

\subsection{Short loxodromic elements}

Let $X$ be a $\delta$-hyperbolic geodesic metric space, and $S$ a finite (not necessarily symmetric) set of isometries of $X$. Denote    $S^{\le n}:=\bigcup_{i=0}^n S^i$, where $S^i$ is the set of isometries that can be expressed as $i$ products of isometries in $S$. Note that $S^{\le n}$ may not be symmetric as well.

The following result improves \cite[Lemma 3.2]{CJY} and allows $S$ to be a finite non-symmetric set of isometries. The proof is almost word by word copy of \cite[Lemma 3.2]{CJY}, but the latter assumes $S$ to be symmetric. So we include a detailed proof here to spell out the differences.
\begin{lemma}\label{ShortHypLem}
Assume that   $\lambda(S,X) \geq 30\delta$. Then $S^{\le 2}$ contains a loxodromic element $b$ with the following property.

Let $o\in X$ such that $|\lambda(S,o)- \lambda(S,X)|<\delta$. There  \revise{exists a  constant $ c_0>0$} depending only on $\delta$ such that
$$
|o-bo|\ge \lambda(S,X)-10\delta
$$
and the path
$$
\alpha :=\bigcup_{n\in \mathbb Z} b^n[o,bo]
$$ is a   \revise{$c_0$}-quasi-axis for $b$.
\end{lemma}

\begin{proof}

Set $L_0=4\delta$ and then  $\lambda(S,o)>7L_0$. Denote by $S_0$ the (non-empty) set of elements $s\in S$ so that $$|o- so|\ge \lambda(S,o)-2L_0-\delta.$$

Let   $t\in S$ such that $\lambda(S,o)=|o- to|$, and $m\in [o, to]$ so that $|o- m| =L_0$.

The main observation is as follows.
\begin{claim}
There exists an isometry $s\in S_0$ such that  $s$ is either loxodromic with $\langle o, s^2o\rangle_{so}\le L_0$  or satisfies
$$
\max\{\langle t^{-1} o, s o \rangle_o,\langle t o, s^{-1} o \rangle_o\}\le L_0.
$$
\end{claim}
\begin{proof}[Proof of the Claim]
Assume to the contrary that for all $s\in S_0$, we have
\begin{equation}\label{MaxiumEQ}
\max\{\langle t^{-1} o, s o \rangle_o,\langle t o, s^{-1} o \rangle_o \}> L_0.
\end{equation}
Moreover, each $s\in S_0$ is either non-loxodromic or loxodromic   with $\langle o, s^2o\rangle_{so}> L_0$.
If $s\in S_0$ is non-loxodromic, by Lemma \ref{LoxoCriterionLem}, we have $$\langle o, s^2o\rangle_{so}\ge |o-so|/2-3\delta\ge (\lambda(S,o)-2L_0-\delta)/2-3\delta\ge L_0.$$ Hence, for each $s\in S_0$, we have $\langle o, s^2o\rangle_{so}\ge L_0$. In particular, $\langle t^{-1}o, to\rangle_{o}\ge L_0$.

\begin{figure}[ht]

\begin{tikzpicture}[x=0.75pt,y=0.75pt,yscale=-1,xscale=1]

\draw    (42,46) -- (320,37) ;
\draw [shift={(320,37)}, rotate = 358.15] [color={rgb, 255:red, 0; green, 0; blue, 0 }  ][fill={rgb, 255:red, 0; green, 0; blue, 0 }  ][line width=0.75]      (0, 0) circle [x radius= 3.35, y radius= 3.35]   ;
\draw [shift={(42,46)}, rotate = 358.15] [color={rgb, 255:red, 0; green, 0; blue, 0 }  ][fill={rgb, 255:red, 0; green, 0; blue, 0 }  ][line width=0.75]      (0, 0) circle [x radius= 3.35, y radius= 3.35]   ;
\draw    (42,46) -- (130,101) -- (226,160) ;
\draw [shift={(226,160)}, rotate = 31.57] [color={rgb, 255:red, 0; green, 0; blue, 0 }  ][fill={rgb, 255:red, 0; green, 0; blue, 0 }  ][line width=0.75]      (0, 0) circle [x radius= 3.35, y radius= 3.35]   ;
\draw [shift={(86,73.5)}, rotate = 32.01] [color={rgb, 255:red, 0; green, 0; blue, 0 }  ][fill={rgb, 255:red, 0; green, 0; blue, 0 }  ][line width=0.75]      (0, 0) circle [x radius= 3.35, y radius= 3.35]   ;
\draw [shift={(178,130.5)}, rotate = 31.57] [color={rgb, 255:red, 0; green, 0; blue, 0 }  ][fill={rgb, 255:red, 0; green, 0; blue, 0 }  ][line width=0.75]      (0, 0) circle [x radius= 3.35, y radius= 3.35]   ;
\draw    (86,73.5) -- (98,45) ;
\draw [shift={(98,45)}, rotate = 292.83] [color={rgb, 255:red, 0; green, 0; blue, 0 }  ][fill={rgb, 255:red, 0; green, 0; blue, 0 }  ][line width=0.75]      (0, 0) circle [x radius= 3.35, y radius= 3.35]   ;
\draw    (178,130.5) -- (171,175) ;
\draw [shift={(171,175)}, rotate = 98.94] [color={rgb, 255:red, 0; green, 0; blue, 0 }  ][fill={rgb, 255:red, 0; green, 0; blue, 0 }  ][line width=0.75]      (0, 0) circle [x radius= 3.35, y radius= 3.35]   ;
\draw    (15,221) -- (226,160) ;
\draw [shift={(226,160)}, rotate = 343.88] [color={rgb, 255:red, 0; green, 0; blue, 0 }  ][fill={rgb, 255:red, 0; green, 0; blue, 0 }  ][line width=0.75]      (0, 0) circle [x radius= 3.35, y radius= 3.35]   ;
\draw [shift={(15,221)}, rotate = 343.88] [color={rgb, 255:red, 0; green, 0; blue, 0 }  ][fill={rgb, 255:red, 0; green, 0; blue, 0 }  ][line width=0.75]      (0, 0) circle [x radius= 3.35, y radius= 3.35]   ;
\draw  [line width=0.75]  (227,157) .. controls (229.53,153.08) and (228.83,149.85) .. (224.91,147.32) -- (216.1,141.65) .. controls (210.5,138.04) and (208.96,134.27) .. (211.49,130.35) .. controls (208.96,134.27) and (204.9,134.43) .. (199.29,130.82)(201.81,132.44) -- (191.68,125.91) .. controls (187.75,123.38) and (184.53,124.08) .. (182,128) ;
\draw  [line width=0.75]  (98,40) .. controls (97.84,35.33) and (95.43,33.08) .. (90.76,33.24) -- (77.88,33.69) .. controls (71.22,33.92) and (67.81,31.7) .. (67.65,27.04) .. controls (67.81,31.7) and (64.56,34.15) .. (57.89,34.38)(60.89,34.27) -- (46.76,34.76) .. controls (42.09,34.92) and (39.84,37.33) .. (40,42) ;
\draw  [line width=0.75]  (38,48) .. controls (35.6,52) and (36.4,55.2) .. (40.4,57.6) -- (48.9,62.7) .. controls (54.61,66.13) and (56.27,69.85) .. (53.87,73.85) .. controls (56.27,69.85) and (60.33,69.56) .. (66.05,72.99)(63.48,71.45) -- (73.4,77.4) .. controls (77.4,79.8) and (80.6,79) .. (83,75) ;

\draw (22,32.4) node [anchor=north west][inner sep=0.75pt]    {$o$};
\draw (309,44.4) node [anchor=north west][inner sep=0.75pt]    {$to$};
\draw (238,148.4) node [anchor=north west][inner sep=0.75pt]    {$s^{-1} o$};
\draw (75,81.4) node [anchor=north west][inner sep=0.75pt]    {$m_{1}$};
\draw (168,110.4) node [anchor=north west][inner sep=0.75pt]    {$m_{2}$};
\draw (104,25.4) node [anchor=north west][inner sep=0.75pt]    {$m$};
\draw (157,181.4) node [anchor=north west][inner sep=0.75pt]    {$s^{-1} m_{1}$};
\draw (9,192.4) node [anchor=north west][inner sep=0.75pt]    {$s^{-2} o$};
\draw (217.43,110.8) node [anchor=north west][inner sep=0.75pt]  [rotate=-30.44]  {$L_{0}$};
\draw (61.08,8.25) node [anchor=north west][inner sep=0.75pt]  [rotate=-355.21]  {$L_{0}$};
\draw (42.96,71.65) node [anchor=north west][inner sep=0.75pt]  [rotate=-32.9]  {$L_{0}$};
\draw (98.5,48.91) node [anchor=north west][inner sep=0.75pt]  [rotate=-23.15]  {$\leq 3\delta $};
\draw (144.68,139.7) node [anchor=north west][inner sep=0.75pt]  [rotate=-3.25]  {$\delta \geq $};

\end{tikzpicture}
\caption{The configuration with $\epsilon=-1$. Replacing $t$,  $s^{-1}$ and $s^{-2}$ by $t^{-1}$, $s$ and $s^2$ respectively gives the other possibility $\epsilon=1$.}
    \label{fig:shortloxodromic}
\end{figure}
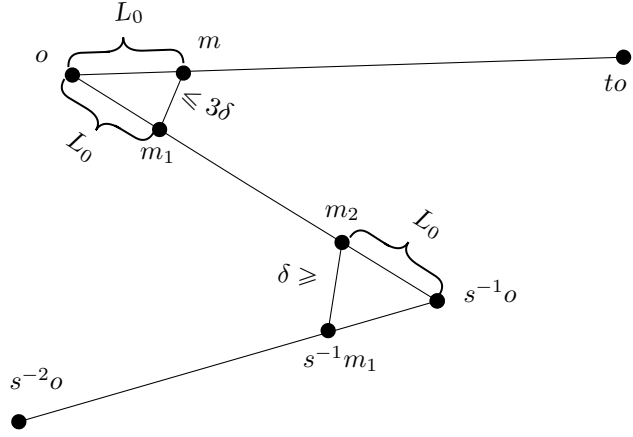
By (\ref{MaxiumEQ}), there exists   $\epsilon\in \{1, -1\}$ so that $\langle t^{-\epsilon} o, s^{\epsilon} o \rangle_o>L_0$. See Figure \ref{fig:shortloxodromic} for illustration with $\epsilon=-1$. Choose $m_1, m_2\in [o, s^\epsilon o]$ for $s\in S_0$ so that $|o- m_1|=|s^\epsilon o- m_2|= L_0$. Using hyperbolicity, \revise{we obtain the case $\epsilon=1$ of   the following inequality}  $$\langle s^\epsilon o, to\rangle_o\ge \min\{\langle s^\epsilon o, t^{-\epsilon} o\rangle_o, \langle t^{-1}o, to\rangle_{o}\}-\delta \ge L_0-\delta$$ \revise{where the other case $\epsilon=-1$ follows from $\langle s^\epsilon o, to\rangle_o=\langle s^\epsilon o, t^{-\epsilon} o\rangle_o$}.   The  $\delta$-thin triangle property then implies  $|m- m_1|\le 3\delta$. Noting $\langle o, s^{2\epsilon} o\rangle_{s^\epsilon o} = \langle s^{-\epsilon} o, s^{\epsilon} o\rangle_{ o} \ge L_0$, we obtain again from    $\delta$-thin triangle  that $|m_2- s^\epsilon m_1|\le \delta$.

We shall derive  $\lambda(S,m)<\lambda(S,X)$, which is a  contradiction. Indeed, for each $s\in S_0,$
\begin{align*}
|m- s^\epsilon m|&\le  2 |m-m_1|+|m_1- s^\epsilon m_1|\le 7\delta+|m_1-m_2|\\
&\le 7\delta+|o-s^\epsilon o|-2L_0 \\
& \le \lambda(S,o)-2L_0+7\delta \le \lambda(S,o)-\delta.
\end{align*}
If $s\in S\setminus S_0$, then $|o- so|\le \lambda(S,o)-2L_0-\delta$ and thus  $$|m- sm|\le 2|o- m|+|o- so|\le 2L_0+|o-so|< \lambda(S,o)-\delta.$$
We obtain the contradiction $\lambda(S,m)\le \lambda(S,o)-\delta< \lambda(S,X)$. The proof of the claim is now complete.
\end{proof}

By the above claim,  there exists $s\in S_0$ such that either
$$\langle s o, s^{-1} o \rangle_o+\delta \le L_0+\delta \le \frac{1}{4} |o-so|$$
or
$$
\max\{\langle t^{-1} o, s o \rangle_o, \langle t  o, s^{-1} o \rangle_o\} +\delta \le L_0+\delta \le \frac{1}{4} \min\{|o-so|, |o- to|\}.
$$
\revise{In the second case, it is possible that $s=t$, say if $S_0=\{t\}$. If this happens, the second inequality is exactly the first one.} Accordingly, we apply Lemma \ref{LoxoCriterion2Lem} to $g=h=s$ or $g=t^{-1}, h=s$. In the first case, $|o-so|>\lambda(S,o)-2L_0-\delta\ge \lambda(S,X)-10\delta$. In the second case, as $\langle t^{-1}o,so\rangle_o\le L_0$, it follows that $|o-tso|=|t^{-1}o-so|\ge |o-t^{-1}o|+|o-so|-2L_0\ge 2\lambda(S,o)-2L_0-\delta\ge \lambda(S,X)-10\delta$. Then $b:=s$ or $b:=ts$ is the desired loxodromic element with the $c_0$-quasi-axis $\alpha$.  The proof is then completed by   Lemma \ref{LoxoCriterion2Lem}.
\end{proof}

If the space $X$ is a tree, the lower bound on the displacement  $\lambda(S,X)$ could be improved as follows, which is equivalent to say that $S$ has no global fixed point.

\begin{corollary}\label{ShortHypLemOnTree}
Let $S$ be a finite set of isometries acting on a simplicial tree $X$ with combinatorial metric.
If      $\lambda(S,X) >0$, then $S^{\le 2}$ contains a loxodromic element $b$ with the following property.

Let $o\in X$ such that $\lambda(S,o)= \lambda(S,X)$. There exists a \revise{uniform} constant $ c_0>0$  such that
$$
|o-bo|\ge \lambda(S,X)-10
$$
and the path
$$
\alpha :=\bigcup_{n\in \mathbb Z} b^n[o,bo]
$$ is a $c_0$-quasi-axis for $b$.
\end{corollary}
\begin{proof}
\revise{As a tree is $0$-hyperbolic, letting $\delta=\min\{\lambda(S,X)/30, 1\}$ in the proof of Lemma \ref{ShortHypLem}   proves the conclusion.}
\end{proof}


\subsection{Proper product actions}

Recall that an isometric group action on a metric space $G\curvearrowright X$ is called \textit{metrically proper} if for any bounded subset $K\subset X$, the set $$\{g\in G\mid gK\cap K\neq \emptyset \}$$ is finite.   As we always consider metrically proper actions in this paper, we shall write proper actions for brevity.

\begin{proposition}\cite[Theorem 6.1]{But22}\label{FiniteIndexProp}
Suppose that $X=\prod_{i=1}^lX_i$ is a finite product of connected graphs with the $\ell_1$-metric or $\ell_\infty$-metric, where each $X_i$ is equipped with the induced path metric. Suppose that $G$ is any group of isometries on $X$. Then $G$ has a finite index subgroup    preserving each factors of $X$.
\end{proposition}




In what follows, we say a \textit{proper product action} of $G$ on $X=\prod_{i=1}^l X_i$ if it acts properly on $X$ and preserves each factors.






\section{Short Loxodromic Element for product actions}\label{ShortLoxoSection}

In this section, we always assume that $G$ admits a proper product action on  \revise{a finite product of} hyperbolic spaces $X=\prod_{i=1}^lX_i$ with either coarsely dense orbits or shadowing property on each factor. Fix a basepoint $o=(o_1,\cdots,o_l)\in X$. Let  $D$ be the coarsely density constant (i.e.: for any $x\in X$, there exists $g\in G$ such that $|x-go|\le D$) or shadowing constant given by Lemma \ref{shadow}. In the latter case, $G$ acts coboundedly on each factor by Lemma \ref{ShadowCoboundedLem}. So we refer the readers that $D$ is independent of $o$ in both cases.

Let $S$ be a finite subset in $G$. For $1\le i\le l$, recall that $\lambda(S,x_i)$ is the displacement of $S$ at $x_i\in X_i$, and $\lambda(S, X_i)$ the displacement of $S$ on $X_i$ as in (\ref{displacement}).   In particular, $$\lambda(S,Go)=\inf_{g\in G}\max_{s\in S}|go-sgo|$$ Note that $\lambda(S, X), \lambda(S,X_i)$ are invariant under conjugation.

Fix once for all a finite generating set $F$ of $G$. Let $\Phi: G\to X$ be the orbital map given by $g\mapsto go$ for any $g\in G$. \revise{Let $[o,go]$ be a fixed choice of geodesic segment between $o$ and $go$.} For each edge $[1,s]$ for $s\in F$ \revise{in $\mathcal G(G,F)$}, define $$
\Phi([1,s]):=[o, so]
$$
which then extends equivariantly to all the edges of Cayley graph $\mathcal G(G,F)$  to $X$. 

{Moreover,}   $$\Phi_i:=\pi_i\cdot \Phi: \;\mathcal G(G,F) \to X_i$$ is induced by the corresponding orbital map of $G$ on $X_i$ with respect to the basepoint $o_i\in X_i$.

The main result of this section is as follows.

\begin{theorem}\label{KeyThm}
For any $M>0$, there exists $N_0=N_0(M, D, \delta, l)>0$ with the following property. For any finite subset $S$ of $G$ with $|S|\geq N_0$, there exists $1\le i \le l$ such that $\lambda(S,X_i) > M$. 
\end{theorem}
We emphasize that $S$ is not required to be symmetric.


The following proposition plays a very important role in proving Theorem \ref{KeyThm}.

\begin{proposition}\label{keyprop}
There exists a constant $C_1=C_1(D,\delta,l)>0$
such that for any $M>0$ and any finite subset $S$ in $G$, one of the following must hold:
\begin{enumerate}
    \item there exists $1\le i\le n$ such that $\lambda(S,X_i)> M.$
    \item $\lambda(S,Go)\le C_1(M+1)$ for any $o\in X$.
\end{enumerate}
\end{proposition}

We first need  the following technical lemma. 
\begin{lemma}\label{keylemma}
Let $S$ be a finite set of isometries  on a hyperbolic space $X$. Fix two points $x, o\in X$.  Then there exists a point $z\in [o,so]$ for some $s\in S$    such that $\lambda(S,z)\leq 6\lambda(S,x)+24\delta$.
\end{lemma}
\begin{proof}
Denote $L=\lambda(S,x)$ in the proof. We first prove the following claim.
\begin{claim}
    For any $s\in S$, if $m\in [o,so]$ is the middle point and  $y\in [o,x]$, $y'\in [so,x]$ are the points satisfying
\begin{align}\label{choiceyEQ}
|x-y| =|x-y'|=  \max\{|x-o| -|o-m|-L/2, 0\},
\end{align}
then $|sy-y|\le L+2\delta$ and $|y-m|\le 5L/2+10\delta$.
\end{claim}
\begin{proof}[Proof of Claim]
We prove in two cases according to the selection of the maximum value in (\ref{choiceyEQ}).

\textbf{Case I}: $|x-o| -|o-m|-L/2\le 0$.
In this case, $x=y=y'$. Clearly $|sx-x|\le L<L+2\delta$. Then we estimate $|x-m|$.

The Gromov product $\langle o,so\rangle_x$ with respect to $x$ is bounded as follows:
\begin{align*}
\langle o,so\rangle_x &=\frac{|x-o|+|x-so|-|o-so|}{2}\\
&\le \frac{2|x-o|+|x-sx|-2|o-m|}{2}\\
&\le |x-o|-|o-m|+L/2\le L,
\end{align*}
from which we then obtain $d(x,[o,so])\le \langle o,so\rangle_x+4\delta\le L+4\delta$. Let $z\in [o,so]$ satisfy $|x-z|=d(x,[o,so])$, then $$\big| |o-x|-|o-z|\big|\le |x-z|\le L+4\delta.$$

As $|o-so|=2|o-m|\le |o-x|+|x-sx|+|sx-so|\le 2|o-x|+L$, one has also
$$\big||o-x|-|o-m|\big| \le L/2.$$

Hence, $|z-m|=\big||o-z|-|o-m|\big|\le 3L/2+4\delta$ and then $|x-m|\le |x-z|+|z-m|\le 5L/2+8\delta$.

\textbf{Case II}: $|x-o| -|o-m|-L/2> 0$.
    At first,  compute the Gromov product of $\langle o, so\rangle_x$ with respect to $x$:
\begin{align*}
\langle o,so\rangle_x &=\frac{|x-o|+|x-so|-|o-so|}{2}\\
&\ge \frac{2|x-o|-|x-sx|-2|o-m|}{2}\\
&\ge |x-o|-|o-m|-L/2.
\end{align*}
From (\ref{choiceyEQ}), $|x-y| =|x-y'|\le \langle o,so\rangle_x$, so by $\delta$-thin triangle property, we have $|y-y'|\le \delta$.

Look at the triangle $\Delta(x,sx,so)$ where the side $[x,sx]$ is of length at most  $L=\lambda(S,x)$. As $\big||so-sy|-|so-y'|\big|= \big||so-sx|-|so-x|\big|\le |x-sx|\le L$,  the $\delta$-thin property implies   $|y'-sy|\le L+\delta$, and then  $|y-sy|\le L+2\delta$. Noting from (\ref{choiceyEQ}) that     $|y-o|=|x-o|-|x-y|=L/2+|o-m|$, it follows that
\begin{align*}
 |o-so|=2|o-m| &\le  |o-y|+|y-sy|+|sy-so|\\
 &\le 2|o-y|+L+2\delta\le  2|o-m|+2L+ 2\delta \\
 &\le |o-so|+2L+2\delta.
\end{align*}
Thus, the path $[o,y][y,sy][sy,so]$ is a $(1,2L+2\delta)$--quasi-geodesic. This implies $\langle o,so\rangle_y\le L+\delta$, so the distance   $d(y,  [o,so])\le \langle o,so\rangle_y +4\delta \le L+5\delta$. As  $\big||o-m|-|o-y|\big|=L/2$ we obtain $$|y-m|\le 2d(y,  [o,so])+L/2\le 5/2L+10\delta.$$
In conclusion, the claim is proved.
\end{proof}

Now, let $t\in S$ be an element maximizing  $\lambda(S,o)$ at the base point $o$:
$$
|to-o|=\max_{s\in S}|so-o|.
$$
If   $\hat m\in [o,to]$ is the middle point, and $\hat y\in [o,x]$ satisfies (\ref{choiceyEQ}) with $m$ replaced with $\hat m$, the claim implies that $|t\hat y-\hat y|\le L+2\delta$ and $|\hat y-\hat m|\le 5/2L+10\delta$.

Recall that $|to-o|$ is maximal among $s\in S$. If $y$ is the point satisfing (\ref{choiceyEQ}) for the element $s\in S\setminus \{t\}$, then  $|x-\hat y|\le |x-y|$, so $\hat y\in [x, y]$. Applying the above claim shows $|y-sy|\le L+2\delta$. By assumption,  for  $|x-sx|\le L$,   the $2\delta$-thin quadrilateral  with four points $x,sx,y,sy$    shows that  $$|\hat y-s\hat y|\le 2\delta+\max\{|x-sx|, |y-sy|\}\le L+4\delta$$ for any $s\in S$. Thus, $$|s\hat m-\hat m|\le |s\hat y-\hat y|+2|\hat y-\hat m|\le 6L+24\delta$$ for all $s\in S$. The proof is complete.
\end{proof}

Now we turn to the proof of Proposition \ref{keyprop}.

\begin{proof}[Proof of Proposition \ref{keyprop}]
For any $M>0$ and any finite set $S$ in $G$, assume  that $\lambda(S, X_i)\leq M$ for all $1\le i\le l$.   We shall show that $\lambda(S, Go)$ is bounded by a linear multiple of $M$ depending only on $l,\delta,D$ for any $o\in X$. Noting that $$\lambda(S, Go)=\lambda(g^{-1}Sg, Go)\leq \lambda(g^{-1}Sg,o)\leq \sum_{i=1}^l\lambda(g^{-1}Sg,o_i),$$
the idea is to find an element $g\in G$ such that $\lambda(g^{-1}Sg,o_i)$ are   uniformly bounded for all $1\le i\le l$.

We construct inductively the conjugates of $S$. Set $S_0=S$ and $M_0=0$ to start.

Suppose $S_i$ has been defined for $i\ge 0$, and      for $1\le i\le l-1$, set
\begin{equation}\label{OldDIsplacement}
M_i := \max\{\lambda(S_i, o_j): 1\le j\le i\}.
\end{equation}
The goal is to inductively define  $S_{i+1}(i\geq 0)$ and then to   bound   $M_{i+1}$ by a linear multiple of $M$ depending on $i, \delta, D$.

As the displacement of $S$ is a conjugacy invariant, we have  $\lambda(S_{i}, X_{i+1})=\lambda(S, X_{i+1})\le M$.  It follows from Lemma \ref{keylemma} that there exists $s_{i}\in S_{i}$ and $z_{i+1}\in [o_{i+1},s_{i}o_{i+1}]$ such that $\lambda(S_{i}, z_{i+1})\leq 6M+24\delta.$

\begin{claim}
    There exists an element $g_{i+1}\in G$, such that $|z_{i+1}-g_{i+1}o_{i+1}|\le D$ and $d(g_{i+1}o_j,[o_j,s_io_j])\le D$ for any $1\le j\le i$.
\end{claim}
\begin{proof}[Proof of Claim]
    We prove the claim in the two cases.

    \textbf{Case I}: $G$ acts on $X$ with coarsely dense orbits. For the chosen $s_i$ and $z_{i+1}$, we choose a point $z$ on the geodesic path $[o,s_io]\subset X$ such that the $(i+1)$-th coordinate of $z$ satisfies $\pi_{i+1}(z)=z_{i+1}$.  As $X\subset N_D(Go)$, there exists $g_{i+1}\in G$ such that $|z-g_{i+1}o|\le D$. Since $|z-g_{i+1}o|=\sum_{1\le j\le l}|\pi_j(z)-g_{i+1}o_j|\le D$, then $|z_{i+1}-g_{i+1}o_{i+1}|\le D$ and $d(g_{i+1}o_j,[o_j,s_io_j])\le d(g_{i+1}o_j,\pi_j(z))\le D$ for any $j\le i$.

    \textbf{Case II}: $G$ acts on each factor with shadowing property.

    Consider a path $\gamma\in \mathcal P$ from $1$ to $s_{i}$ in the Cayley graph $\mathcal{G}(G,F) $. By Lemma \ref{shadow}, $[o_{i+1},s_{i}o_{i+1}]\subset N_D(\Phi_{i+1}(\gamma))$. Hence, there exists $g_{i+1}\in \gamma$ such that $|z_{i+1}-g_{i+1}o_{i+1}|\leq D$.

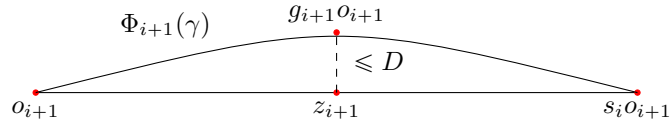
\begin{figure}[ht]
  \centering
  \begin{tikzpicture}
\draw [black] (-4,0)--(4,0);
\filldraw [red] (-4,0) circle (1pt) node[black, below] {$o_{i+1}$};
\filldraw [red] (0,0) circle (1pt) node[black, below] {$z_{i+1}$};
\filldraw [red] (4,0) circle (1pt) node[black, below] {$s_{i}o_{i+1}$};
\filldraw [red] (0,0.8) circle (1pt) node[black, above] {$g_{i+1}o_{i+1}$};
\draw [black] (-4,0) .. controls (0,1) .. (4,0);
\draw [dashed, black] (0,0)--(0,0.8);
\draw (0.1,0.4) node[black,right] {$\leq D$};
\draw (-3,0.9) node[black,right] {$\Phi_{i+1}(\gamma)$};
\end{tikzpicture}
  \caption{The shadowing quasi-geodesic $\Phi_{i+1}(\gamma)$ in  $X_{i+1}$}
\end{figure}

For each $j\le i$, the shadowing property of $G\curvearrowright X_j$ implies that  $\Phi_j(\gamma)$ is a uniform quasi-geodesic connecting $o_j$ to $s_{i}o_j$ in $X_j$, so by Lemma \ref{shadow},  $d(g_{i+1}o_j,[o_j,s_{i}o_j])\leq D$. See Figure \ref{ShadowQG} for illustration.

In conclusion, the claim is proved.
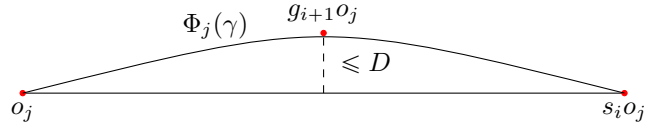
\begin{figure}[ht]
  \centering
  \begin{tikzpicture}
\draw [black] (-4,0)--(4,0);
\filldraw [red] (-4,0) circle (1pt) node[black, below] {$o_j$};
\filldraw [red] (4,0) circle (1pt) node[black, below] {$s_{i}o_j$};
\filldraw [red] (0,0.8) circle (1pt) node[black, above] {$g_{i+1}o_j$};
\draw [black] (-4,0) .. controls (0,1) .. (4,0);
\draw [dashed, black] (0,0)--(0,0.8);
\draw (0.1,0.4) node[black,right] {$\leq D$};
\draw (-2,0.9) node[black,right] {$\Phi_j(\gamma)$};
\end{tikzpicture}
  \caption{The shadowing quasi-geodesic $\Phi_j(\gamma)$ in   $X_j$ for $1\le j\le i$}
\end{figure}\label{ShadowQG}
\end{proof}

Next, we define the new conjugate   $S_{i+1}:=g_{i+1}^{-1}S_{i}g_{i+1}$. We need to estimate $M_{i+1}$ in the definition (\ref{OldDIsplacement}) where $S_i$ in the definition of $M_i$ is changed to $S_{i+1}$.

First of all, we consider the displacement on $o_{i+1}$ in $X_{i+1}$:
\begin{align*}
    \lambda(S_{i+1},o_{i+1}) & =\max_{t\in S_{i+1}}|o_{i+1}-to_{i+1}|=\max_{t\in S_{i}}|g_{i+1}o_{i+1}-tg_{i+1}o_{i+1}|\\
    &\leq 2|g_{i+1}o_{i+1}-z_{i+1}|+ \max_{t\in S_{i}}|z_{i+1}-tz_{i+1}|\\ &\le \lambda(S_{i},z_{i+1})+2D \leq 6M+24\delta + 2D.\\
\end{align*}

It remains  to update the displacement of $S_{i+1}$ on the point $o_j\in X_j$ for each $1\le j\le i$. From Figure \ref{ShadowQG}, $|o_j-g_{i+1}o_j|\le D+|o_j-s_io_j|\le \lambda(S_i,o_j)+D$ as $s_i\in S_i$.
Then we obtain
\begin{align*}
    \lambda(S_{i+1},o_j) & =\max_{t\in S_{i+1}}|o_j-to_j|
    =\max_{t\in S_{i}}|g_{i+1}o_j-tg_{i+1}o_j|\\
    &\leq \max_{t\in S_{i}}|o_j-to_j|+2|o_j-g_{i+1}o_j|\\
    &\leq 3\lambda(S_{i},o_j)+2D \leq 3 M_i +2D.\\
\end{align*}
Thus, we get
\begin{equation}\label{LinearG}
    M_{i+1}\leq \max\{6M+27\delta+2D,3M_i+2D\}.
\end{equation}


After $l$ steps, we obtain $S_{l}=g^{-1}Sg$ where  $g=g_1g_2\cdots g_l$ and      $$\lambda(S, Go)=\lambda(S_{l}, Go)\le \lambda(S_{l},o)\le \sum_{i=1}^l\lambda(S_{l},o_i)\le l M_l.$$
It follows from (\ref{LinearG}) that there exists a constant $C_1=C_1(l,\delta,D)$ such that $lM_l\le C_1(M+1)$. The proof is completed.
\end{proof}

The proof of Proposition \ref{keyprop} actually implies the following lemma. 

\begin{lemma}\label{Equi}
    For any finite $S\subset G$ and $o\in X$, we have $$\lambda(S,X)\le \lambda(S,Go)\le C_1\lambda(S,X)+C_1$$
\end{lemma}
\begin{proof}
    As $Go\subset X$, $\lambda(S,X)\le \lambda(S,Go)$ follows from the definition (\ref{displacement}) of displacement.

     Pick $o'=(o_1',\cdots,o_l')\in X$ such that $\lambda(S,X)=\lambda(S,o')$. For each $1\le i\le l$, let $s_i\in S$ so that $\lambda(S,o_i')=|o_i'-s_io_i'|$. Then the second inequality follows from Proposition \ref{keyprop} combined with the following estimate:
    $$\lambda(S,X_i)\le \lambda(S,o_i')=|o_i'-s_io_i'|\le \sum_{1\le j\le l}|o_j'-s_io_j'|\le \lambda(S,o')=\lambda(S,X).$$

\end{proof}

The following lemma comes from the proper product action $G\curvearrowright X$ and Lemma \ref{Equi}.

\begin{lemma}\label{MN}
For any $M > 0$ there exists $N = N(M) > 0$ such that for any finite subset $S$ of $G$, if $|S| > N$, then $\lambda(S,X) > M.$
\end{lemma}
\begin{proof}
Since $G$ acts   properly on $X$, for any $M>0$ and $o\in X$, the following set $$N(o, C_1M+C_1):=\{g\in G: |o-go|\leq C_1M+C_1\}$$ has  finite cardinality. Let $N:=\inf_{o\in X}|N(o, C_1M+C_1)|=|N(o', C_1M+C_1)|$ for some $o'\in X$.

For any finite subset $S\subset G$ with $|S|>N$, the displacement of $S$ on $o'$ is larger than $C_1M+C_1$. As $|g^{-1}Sg|=|S|>N$, $\lambda(g^{-1}Sg,o')>C_1M+C_1$ for any $g\in G$. Hence, $\lambda(S,Go')=\inf_{g\in G}\lambda(g^{-1}Sg,o')>C_1M+C_1$. By Lemma \ref{Equi}, $\lambda(S,X)\ge \frac{1}{C_1}\lambda(S,Go')-1>M$.
\end{proof}

Now we can complete the proof of Theorem \ref{KeyThm}.
\begin{proof}[Proof of Theorem \ref{KeyThm}]
Let  $C_1$ be the constant  given by   Proposition \ref{keyprop}. For any $M>0$, let $N_0=N(C_1M+C_1)$ be given by Lemma \ref{MN}.    For any finite set $S$  with $|S|\ge N_0$,  Lemma \ref{MN} gives that  $\lambda(S,X)>C_1M+C_1$. Thus, the displacement of $S$ must be larger than $M$ on some factor space $X_i$ by     Proposition \ref{keyprop}. This is what we wanted to prove.
\end{proof}

\section{Uniform Exponential Growth}\label{UEGSection}

The goal of this section is Theorem \ref{Thm1}, re-stated as follows, that groups with proper product actions on hyperbolic spaces have locally uniform exponential growth.
\begin{theorem}\label{ueg}
Suppose that $G$ admits a proper product action  on hyperbolic spaces with either coarsely dense orbits or shadowing property on each factor. Then there exists a constant $\omega_0>1$ such that $\omega(H)\ge \omega_0$ for any finitely generated  non-virtually abelian subgroup $H$.
\end{theorem}


At first, we observe the following elementary facts. For any two subsets $H,K\subset G$, denote $$\mathcal C(H,K):=\{[h,k]: h\in H, k\in K\}$$ and $[H,K]$ is the group generated by $\mathcal C(H,K)$.
\begin{lemma}\label{FiniteCenter}
Let $G$ be a  group generated by a finite   set   $S$. Suppose that $g\in G$ is   of infinite order so that $\langle g\rangle\cap Z(G)$ is finite. Then  there exists $a\in S$ such that $[g^k,a]\neq 1$ and $[g^k,a]\ne [g^j,a]$  for any two distinct integers $k \ne j \geq 1$. In particular,    $\mathcal C(G,G)$ is infinite.
\end{lemma}
\begin{proof}
Assume to the contrary that each element of $S$ commutes with some nontrivial power of $g$. Let $K>0$ be the least common multiplier so that $g^K$ commutes with every element in $S$. Thus, $\langle g^K\rangle\leq Z(G)$, which contradicts the finiteness of $\langle g\rangle \cap Z(G)$. Hence, there exists $a\in S$ such that $[g^k,a]\neq 1$ for all   $k\geq 1$. This implies  $[g^i,a]\neq [g^j,a]$ for any $i\ne j$.  Thus, $\mathcal C(G,G)$ is infinite.
\end{proof}

\revise{The following lemma appears to be  a well-known fact, though we could not locate a precise reference. For completeness, we include a simple proof here.

\begin{lemma}\label{NormalGene}
    Let $G$ be a group with a generating set $S$. Then $[G,G]$ is normally generated by $\mathcal C(S,S)$. That is, $[G,G]$ is the normal closure of $\mathcal C(S,S)$.
\end{lemma}
\begin{proof}
    For any $g,h\in G$, denote $g^h:=hgh^{-1}$. Then for any $g,h,f\in G$, one has that
    $$[g,hf]=ghfg^{-1}(hf)^{-1}=ghg^{-1}h^{-1}hgfg^{-1}f^{-1}h^{-1}=[g,h][g,f]^h$$
    and
    $$[gh,f]=ghf(gh)^{-1}f^{-1}=ghfh^{-1}f^{-1}g^{-1}gfg^{-1}f^{-1}=[h,f]^g[g,f].$$
    Hence, for any $g,h\in G$ with $g=s_1\cdots s_n$ and $h=s_1'\cdots s_m'$ where $s_i,s_j'\in S, 1\le i\le n, 1\le j\le m$, $[g,h]$ can be written as a product of conjugates of $[s_i,s_j']$ for $1\le i\le n,1\le j\le m$.
    As $[G,G]$ is generated by all commutators in $G$, the conclusion follows.
\end{proof}}

\revise{Recall that the \textit{lower central series} of a group $G$ is defined as
$$G_0=G, G_1=[G,G_0], G_2=[G,G_1], \ldots, G_n=[G,G_{n-1}],\ldots.$$
The group $G$ is called \textit{nilpotent} of class $n\ge 1$ if $G_n=\{1\}$ but $G_{n-1}\neq \{1\}$.}

\begin{lemma}\label{LargeCommutator}
Let $G$ be a non-virtually abelian group generated by a finite set $S$. Then there exists a constant $c\in (0,1)$ independent of $S$ such that $|\mathcal C(S^{\le n}, S^{\le n})|\ge cn$ for any $n\in \mathbb N$.
\end{lemma}
\begin{proof}
For simplicity, denote $\mathcal C_n=\mathcal C(S^{\le n}, S^{\le n})$ in this proof. Since $G$ is \revise{not virtually abelian},  we see that $|\mathcal C_1|\ge 2$ and $\mathcal C_n \subset \mathcal C_{n+1}$ for any $n\in\mathbb N$. According to the nilpotent class of $G$, we divide the proof into three cases.

\textbf{Case I}: $G$ is a nilpotent group of class $\le 2$, i.e. $[G,[G,G]]=1$.

\begin{claim}
    $|\mathcal C_n|\ge n$ for any $n\in \mathbb N$.
\end{claim}
\begin{proof}[Proof of Claim]
As $G$ is not virtually abelian and the center $Z(G)$ is abelian, we see that $G/Z(G)$ must be  infinite. By assumption, $[G,G]$ is contained in $Z(G)$, so  $G/Z(G)$ is an infinite abelian group.

Let $\pi: G\to G/Z(G)$ be the quotient map. \revise{Since an infinite abelian group can not be generated by a finite set of torsion elements, the generating set $\pi(S)$ of $G/Z(G)$ contains an infinite order element $\overline a$.} For concreteness, let $a\in S$ so that $\pi(a)=\overline a$. Observe that $\langle a\rangle \cap Z(G)=\{1\}$: \revise{indeed, if   $\langle a\rangle \cap Z(G)$ is of finite index in $\langle a\rangle $, then $\langle \bar a\rangle=\langle a\rangle / Z(G)$ is a finite cyclic group}. By Lemma \ref{FiniteCenter}, there exists $b\in S$ such that $|\{[a^k,b]: 1\le k\le n\}|=n$.

Hence, $|\mathcal C_n|\ge |\{[a^k,b]: 1\le k\le n\}|=n$ and the claim follows.
\end{proof}

\textbf{Case II}: $G$ has a subgroup $H$ of finite index $d$ which is nilpotent of class $\le 2$.

\begin{claim}
    $|\mathcal C_n|\ge \frac{1}{3d}n$ for any $n\in \mathbb N$.
\end{claim}
\begin{proof}[Proof of Claim]
Denote $U:=S^{\le 2d-1}\cap H$. By Shalen-Wagreich Proposition \cite[Proposition 3.3]{ShW}, $U$ generates $H$. As a consequence of Case I, $|\mathcal C(U^{\le n}, U^{\le n})|\ge n$.

Then, if $n\ge 3(2d-1)$, one has
$$|\mathcal C_n|\ge |\mathcal C(U^{\le n/(2d-1)}, U^{\le n/(2d-1)})|\ge \frac{n}{2d-1}-1\ge \frac{n}{3d}.$$
Otherwise, for $1\le n< 3(2d-1)$, we have $|\mathcal C_n|\ge 2\ge \frac{n}{3d}$ as well. Hence,  the claim is proved.
\end{proof}
\textbf{Case III}: $G$ does  not contain a finite index  nilpotent group of class $\le 2$.

\begin{claim}
   $|\mathcal C_n|\ge \frac{1}{2}n$ for any $n\in \mathbb N$.
\end{claim}
\begin{proof}[Proof of Claim]
We argue by contradiction. As the claim clearly holds for $n\le 2$, we suppose that $|\mathcal C_n|< \frac{n}{2}$ for some $n\ge 3$. Note   the following   relation
$$1=|\mathcal C_0|\le |\mathcal C_1|\le \cdots \le |\mathcal C_n|< \frac{n}{2}.$$
By pigeonhole principle, there exists $1\le m\le n-2$ such that $|\mathcal C_m|=|\mathcal C_{m+2}|$. As $\mathcal C_m\subset \mathcal C_{m+2}$, we obtain $\mathcal C_m=\mathcal C_{m+2}$.

For any $s\in S$ and $g,h\in S^{\le m}$, we note that $$s[g,h]s^{-1}=[sgs^{-1},shs^{-1}]\in \mathcal C_{m+2}=\mathcal C_m$$
implying that $\mathcal C_m$ is invariant under conjugacy. \revise{By Lemma \ref{NormalGene}, $[G,G]$ is normally generated by $\mathcal C_1=\{ghg^{-1}h^{-1}: g,h\in S\}$. As  $\mathcal C_m\supset \mathcal C_1$ is conjugacy-invariant, $\mathcal C_m$ generates the commutator subgroup $[G,G]$. }

Let   $H$ be the kernel of the group action of $G$ on $\mathcal C_m$ by conjugation. That is to say, every $h\in H$ commutes all elements in $\mathcal C_m$, so \revise{  $[H,\langle \mathcal C_m\rangle ]=[H,[G,G]]=1$ follows.} In particular, $[H,[H,H]]=1$  shows that $H$ is a nilpotent group of class $\le 2$. This results in a contradiction, as $H$ is of finite index in $G$. Hence, the claim is proved.
\end{proof}
To conclude,  setting $c=\min\{{1}/{3d},{1}/{2}\}$  completes the proof.
\end{proof}

The following lemma is a direct consequence of Lemma \ref{LargeCommutator}, generalizing a result \cite[Lemma II.7.9]{BH} that if a group has finite commutator subgroup, then it is  virtually abelian. Recall that $\mathcal C(G,G)$ denotes the set of all commutators over $G$, so generates $[G,G]$. 
\begin{lemma}\label{FiniteCommutator}
Let $G$ be a finitely generated group so that $\mathcal C(G,G)$ is a finite set. Then $G$ is virtually abelian.
\end{lemma}


As a corollary, we have that

\begin{corollary}\label{AllLineal}
Suppose that a finitely generated infinite group $G$   \revise{admits a proper product action on a finite product of quasi-lines such that $G$ acts on each factor by orientable lineal actions.} Then $\mathcal C(G,G)$ is a finite set and $G$ is virtually abelian with infinite center.
\end{corollary}
\begin{proof}
By Lemma \ref{bounded}, the orbit of the subset $\mathcal C(G,G)$ has finite diameter on every factor. Thanks to the proper action, $\mathcal C(G,G)$ is finite. By Lemma \ref{FiniteCommutator}, $G$ is virtually abelian. As any finitely generated abelian torsion group is a finite group, $G$ must contain an infinite order element.  As a result of Lemma \ref{FiniteCenter},   $Z(G)$ is infinite, otherwise we would obtain that  $\mathcal C(G,G)$ is infinite. The proof is completed.
\end{proof}

Now, we   \revise{strengthen} the conclusion of  Theorem \ref{KeyThm}.

\begin{proposition}\label{refine}
Let $G \curvearrowright \prod_{i=1}^lX_i$ be a proper product action on \revise{a product of } hyperbolic spaces with either coarsely dense orbits or shadowing property on each factor. Assume that $H$ is a non-virtually abelian subgroup of $G$ so that all lineal actions of $H$ on factors are    orientable.   Then for any sufficiently large $M>0$, there exists $N=N(M, D, \delta, l)>0$ \revise{where $D,\delta$ are the same constants as in Theorem \ref{KeyThm}} with the following property.

For any finite  generating set $S$ of $H$, there exists $1\le i \le  l$ such that $$\lambda(\mathcal C(S^{\leq N},S^{\leq N}),X_i) > M$$ and   $H\curvearrowright X_i$ is either focal or of general type.
\end{proposition}

\begin{proof}
Fix any $M\ge \max\{C,30\delta\}$, where $C$ is the constant in Lemma \ref{bounded}.  Let $N_0=N_0(M,D,\delta,l)$ be given by Theorem \ref{KeyThm} and $c$ be the constant in Lemma \ref{LargeCommutator}. \revise{For any $r\ge 0$, denote by $[r]$ the largest integer that does not exceed $r$.} Let $N:=[N_0/c]+1$.

Set $S':=\mathcal C(S^{\leq N},S^{\leq N})$. Then Lemma \ref{LargeCommutator} gives that $|S'|\geq N_0$.  As a result of Theorem \ref{KeyThm}, there exists $1\le i\le l$ such that $\lambda(S',X_i)>M$. As $M>30\delta$, Lemma \ref{ShortHypLem} guarantees the existence of a loxodromic element in $S'^{\le 2}$ on $X_i$. Thus, the action of $H$ on  $X_i$ is either lineal and orientable,  or focal, or of general type.

By way of contradiction,  assume that  $H\curvearrowright X_i$ is lineal and orientable. By the above definition,  $S'$  consists of commutators in $G$, so according to Lemma \ref{bounded}, there exists $x_i\in X_i$ such that $|sx_i-x_i|\leq C$ for any $s\in S'$. This results in   a contradiction as follows
\begin{equation*}
   M<\lambda(S',X_i)=\inf_{x\in X_i}\max_{s\in S'}|sx-x|\leq \max_{s\in S'}|sx_i-x_i|\leq C.
\end{equation*}
We thus prove that $H\curvearrowright X_i$ must be focal or of general type.
\end{proof}

The last ingredient is the following criterion to generate a  free semigroup. Recall that $\tau(g)$ is the stable translation length of any loxodromic isometry $g$.

\begin{proposition}\cite[Prop. 11.1]{BF18}\label{Ping-Pong}
    Let $X$ be a $\delta$-hyperbolic space for some  $\delta\ge 0$. Suppose $g$ and $h$ are loxodromic isometries whose fixed point sets in $\partial X$ are not equal (possibly with a common point). If   $\tau(g),\tau(h)> 10000\delta$, then some pair in $\{g^{\pm 1}, h^{\pm 1}\}$ generates a free semigroup.
\end{proposition}

We are ready to prove Theorem \ref{ueg}.
\begin{proof}[Proof of Theorem \ref{ueg}]
Suppose $G$ admits a proper product action on finite hyperbolic spaces $\prod_{i=1}^l X_i$ with either coarsely dense orbits or shadowing property on each factor.
Let  $H$ be any finitely generated non-virtually abelian subgroup. Choose a subgroup $H_0$ of finite index in $H$ such that all lineal actions of $H_0$ on factors $X_i$  are orientable. As $H_0$ is still non-virtually abelian, it follows from Corollary \ref{AllLineal} that not all $H_0$-actions on factors are lineal and thus $[H:H_0]\le 2^{l-1}$. As a result of  \cite[Prop. 2.4]{H} which gives the relation of uniform exponential growth between a group and its finite index subgroup, $\omega(H)^{2^l-1}\ge \omega(H_0)$.

Fix $N_1$ in Proposition \ref{refine} for $M=10026\delta+C$ where $C$ is the constant in Lemma \ref{bounded}, and then set $N_2=4N_1$. Let $S$ be a  finite symmetric generating set of $H_0$. Thus,  $\mathcal C(S^{\le N_1}, S^{\le N_1})\subset S^{\leq N_2}$. As a result of Proposition \ref{refine}, there exists $1\le i\le l$ such that $\lambda(S^{\leq N_2},X_i)>10026\delta+C$ and $H_0\curvearrowright X_i$ is focal or of general type.

It follows from Lemma \ref{ShortHypLem} that there exists a loxodromic element $g\in S^{\leq 2N_2}$ on $X_i$. Moreover, Lemma \ref{ShortHypLem} and \cite[Lemma 2.7]{CJY} give a lower bound of $\lambda(g,X_i)$ by $\lambda(S^{\leq N_2},X_i)-10\delta$. Hence, $\lambda(g,X_i)>10016\delta+C$ and then $\tau(g)> 10000\delta+C$ by Lemma \ref{StableTransLength}.

Since the action $H_0\curvearrowright X_i$ is not lineal, there exists at least one $h\in S$ such that $h\notin E(g)$. Note that $\tau(h^{-1}gh)=\tau(g)>10000\delta+C$.
By applying Proposition \ref{Ping-Pong}, one gets that some pair in $\{ g^{\pm 1}, h^{-1}g^{\pm 1}h\}$ generates a free semigroup.

To conclude, set $\kappa:=2N_2+2$. Then
\begin{equation*}
    \omega(H_0,S)\geq \omega(H_0,S^{\leq \kappa})^{1/\kappa}\geq \omega(H_0,\{g^{\pm 1}, h^{-1}g^{\pm 1}h\})^{1/\kappa}= 2^{1/\kappa},
\end{equation*}
for any finite symmetric generating set $S$ of $H_0$ and thus $\omega(H_0)\ge 2^{1/\kappa}$. By setting $\omega_0=2^{1/[\kappa (2^l-1)]}$, one has that $\omega(H)\ge \omega(H_0)^{1/(2^l-1)}\ge \omega_0$. The proof is   complete.
\end{proof}


The first application of Theorem \ref{ueg} is given to the class of  hierarchically hyperbolic groups (HHGs).

\begin{corollary}
    Let $(G,\mathfrak G)$ be a  hierarchically hyperbolic group with a BBF-coloring. Then $G$ has locally uniform exponential growth.
\end{corollary}
\begin{proof}
    A result of M. Hagen and H. Petyt \cite[Theorem 3.1, Proposition 3.10]{HP} shows that any hierarchically hyperbolic group $(G,\mathfrak G)$ with a BBF-coloring $\mathfrak G=\bigsqcup_{i=1}^{l} \mathfrak G_i$ admits a quasimedian quasi-isometric embedding into a finite product of hyperbolic spaces $X=\prod_{i=1}^{l} \mathcal C \mathfrak G_i$ and preserves each factor with shadowing property. \revise{Note that if a group $G$ acts isometrically on a finite product of metric spaces so that the orbit map is a quasi-isometric embedding, then the product action must be proper. Therefore, the conclusion then follows immediately from Theorem \ref{ueg}.}
\end{proof}

Another application is given to Burger-Mozes-Wise (BMW) groups. Recall that a group $G$ is a BMW group if it acts by isometries on the product of two trees $T_1\times T_2$ such that every element preserves the product decomposition and the action on the vertex set of $T_1\times T_2$ is free and transitive.

\begin{corollary}
    BMW groups have locally uniform exponential growth.
\end{corollary}

\section{Product Set Growth}\label{PSGSection}

In this section, we shall study a  growth property called \textit{product set growth}, stronger than the purely exponential growth, for groups with a proper product action on a finite product of hyperbolic spaces.

By convention,   we say that a finite set $U$ generates a group $H$ as a semi-group (resp. group), if any element of $H$ can be written as a word over the alphabet $U$ (resp. $U\cup U^{-1}$). 


For any finite subset $U\subset H$ and $n\in \mathbb N$, the \textit{n-th product set} of $U$ is defined as follows $$U^n:=\{s_1\cdots s_i\cdots s_n: s_i \in U\}$$
\revise{Using the above notation, if $U$ generates $H$ as a semi-group (resp. group), we have $H=\bigcup_{n\ge 0} U^n$ (resp. $H=\bigcup_{n\ge 0} (U\cup U^{-1})^n$).}


\revise{
\begin{defn}\cite[Definition 2.0.5]{Ker}
    Let $G$ be an infinite finitely generated group. We say that $G$ has \textit{uniform product set growth} if there exist $\alpha, \beta>0$ such that for every finite (symmetric) generating set $U$ of $G$ we have that $$|U^n|\geq (\alpha |U|)^{\beta n}$$ for any $n\ge 1.$
\end{defn}

The goal of this section, however, is not just to prove this property for the group itself, but try to find a dichotomy for its finitely generated subgroups.

\begin{question}\cite[Question 1.0.1]{Ker}
    For a group $G$, does there exist a class of subgroups $\mathcal H$ such that the following hold?
    \begin{itemize}
        \item Every finitely generated subgroup $H\notin \mathcal H$ has uniform product set growth, for the same $\alpha, \beta>0$.
        \item No finitely generated subgroup $H\in \mathcal H$ has uniform product set growth, for any $\alpha,\beta>0$.
    \end{itemize}
\end{question}
}

\begin{remark}\label{increasing}
     For any finite subset $U$ and $n\geq 1$, $|U^{n}|\geq |U^{n-1}|$. This follows immediately from that $|U^n|\geq |aU^{n-1}|=|U^{n-1}|$ for any $a\in U$.
\end{remark}

Now, we re-state our main result, Theorem \ref{Thm2}, from the  Introduction. Recall that for a loxodromic element $g$ on a hyperbolic space, $E(g)$ is  the  subgroup of elements which fix setwise the endpoints of the quasi-axes $L_g$ of $g$.

\begin{theorem}\label{PSG}
Let $G$ admit  a proper product  action on \revise{a finite product of }  hyperbolic spaces $X=\prod_{i=1}^lX_i$  with shadowing property on factors. Suppose that $G$ acts weakly acylindrically on each factor. Then there exist $\alpha,\beta >0$ depending only on $\delta,l$ such that for every finite symmetric subset $U\subset G$, at least one of the following  is true:
\begin{enumerate}
    \item $\langle U\rangle$ is virtually abelian.
    \item There exists a factor $X_i$ so that no uniform bound exists on $|E(h)\cap E(b)\cap \langle U\rangle|$ for any two independent loxodromic elements $b,h\in \langle U\rangle$ on $X_i$.
    \item $|U^n|\geq (\alpha|U|)^{\beta n}$ for every $n\in \mathbb N$.
\end{enumerate}
Moreover, if there is no lineal action of $\langle U\rangle$ on factors, then $U$ can be non-symmetric and generates a group $\langle U\rangle$ as a semi-group.
\end{theorem}

In general,  a group with uniform exponential growth may not have  uniform product set growth. For example, Kerr showed that any finitely generated group with infinite center (e.g. direct product with abelian group) fails to have uniform product set growth \cite[Cor. 2.2.2]{Ker}.

\subsection{A non-symmetric version of Shalen-Wagreich's result}

\revise{In \cite[Proposition 3.3]{ShW},  Shalen-Wagreich proved a very useful fact that if a symmetric set $S$ generates a group $G$ and $H<G$ is a subgroup of finite index $d$, then $S^{\le 2d-1}$ contains a symmetric generating set of $H$.   The goal of this subsection is to present a non-symmetric version of this result; that is, $S$ is not necessarily symmetric. This shall be used in treating the product set growth for a non-symmetric set.}

\revise{We make use of the classical Reidemeister-Schreier method to get the following lemma. For details about Reidemeister-Schreier method, see \cite[Section II.4]{LS}.}

\begin{lemma}\label{SchreierNormal}
    Let $G$ be a group generated by a finite set $U$  as a semi-group. Let $H$ be a normal subgroup of index $d$ in $G$. Then there exists a finite subset $W\subset U^{\le d^2-d+1}$ such that $W$ generates $H$ as a group. Moreover, $|W|\ge \frac{1}{d}|U|$.
\end{lemma}
\begin{proof}
    As $U$ may be non-symmetric, we re-define the Cayley graph $\mathcal G(G,U)$ as a directed graph: the vertex set consists of all elements of $G$ and there is a directed edge from $g$ to $h$ if $h=gs$ for some $s\in U$. As $U$ generates $G$ as a semi-group, for any $g,h\in G$, there is a directed path from $g$ to $h$ in $\mathcal G(G,U)$. Denote by $d_U(g,h)$ the shortest length of   directed paths from $g$ to $h$ and $d_U(A,B):=\inf_{g\in A, h\in B}d_U(g,h)$ \revise{for any two subsets $A,B\subset G$}. Since $[G:H]=d$, there exists a subset $T=\{a_1,\ldots,a_d\}\subset G$ such that $G=\bigsqcup_{i=1}^da_iH$ and $d_U(1,a_i)=d_U(1,a_iH)$ for $1\le i\le d$.  Set $a_1=1$.

    Denote $\pi: G\to \overline G:=G/H$ and $\overline g:=\pi(g)$. \revise{Let $\overline U:=\pi(U)$. So $\overline U$ generates $\overline G$. Observe that the directed distance from $\overline 1$ to $\overline a_i$ in  $\mathcal G(\overline G,\overline U)$ equals $d_U(1, a_iH)$.} As $\overline G$ has $d$ elements,  $d_U(1,a_i)=d_U(1, a_iH)=d_{\overline U}(\bar 1, \bar a_i)\le d-1$.

    Define a map: $\phi: G\to T$ so that \revise{for every element $g$ in the left coset $a_iH$, $\phi(g)=a_i$. Equivalently speaking, } for any $g\in G$,  the image $\phi(g)\in T$  is uniquely determined by the relation   $\phi(g)^{-1}g\in H$.  Set $W_1:=\{\phi(ua)^{-1}ua: u\in U, a\in T\}\subset H$.

    \begin{claim}
        $W_1$ generates $H$ as a semi-group.
    \end{claim}
    \begin{proof}[Proof of the Claim]
        For any $h\in H$ with $d_U(1,h)=n$, denote $h=u_{n}\cdots u_{1}$ where $u_i\in U$ for $1\le i\le n$. \revise{Note that $a_1=1$. Define $a_{k_1}=\phi(u_{1}a_1)\in T$ and recursively $a_{k_i}=\phi(u_{i}a_{k_{i}})\in T$ for $1\le i\le n-1$. We then apply the following sequence of transformations:
        \begin{align*}
            h=& u_{n}\cdots u_{1}a_1\\
            =& u_{n}\cdots u_{2}(a_{k_1}\phi(u_1a_1)^{-1})u_1a_1\\
            =& u_{n}\cdots u_{2}a_{k_1}(\phi(u_1a_1)^{-1}u_1a_1)\\
            =& \cdots\\
            =& a_{k_n}(\phi(u_{n-1}a_{k_{n-1}})^{-1}u_{n-1}a_{k_{n-1}})\cdots (\phi(u_1a_1)^{-1}u_1a_1).
        \end{align*} }
        As $h\in H$ \revise{and each $\phi(u_ia_{k_{i}})^{-1}u_ia_{k_{i}}\in W_1\subset H$, the element $a_{k_n}$ must be in $H$. Since $a_{k_n}\in T$ by definition is a  left $H$-coset representative, we conclude that  $a_{k_n}=a_1=1$. Hence, } any $h\in H$ can be written as a word over $W_1$. The claim is proved.
    \end{proof}

    \revise {Recall that $G/H$ is a finite quotient group of order $d$. Then} $(\overline a)^d=\overline 1\in G/H$, which implies that $a^d\in H$ for each $a\in T$. Thus, $\phi(ua)^{d-1}ua=\phi(ua)^{d}\cdot \phi(ua)^{-1}ua\in H$ for $u\in U, a\in T$.

    Define $W:=\{a^d: a\in T\}\cup\{\phi(ua)^{d-1}ua: u\in U, a\in T\}$. \revise{Recall that the image of $\phi$ lies in $T$ and for every element $a\in T$, one has $d_U(1,a)\le d-1$. Thus, $d_U(1,a^d)\le d(d-1)$ and $d_U(1,\phi(ua)^{d-1}ua)\le d_U(1,\phi(ua)^{d-1})+d_U(1,u)+d_U(1,a)\le d(d-1)+1$ for any $u\in U, a\in T$}. It follows that $W\subset U^{\le d^2-d+1}$. Moreover, $$\phi(ua)^{-1} ua= (\phi(ua))^{-d} \phi(ua)^{d-1}(ua),$$ so    the above Claim implies that $W$ generates $H$ as a group.

    To give a bound on $|W|$,  note first that \revise{$a_1=1\in T$ and then} $W_2:=\{\phi(u)^{d-1}u: u\in U\}\subset W$. Recall $\phi|_U$ maps $U$ to the set $T$ of $d$ elements, so by   the pigeonhole principle, there is a subset $U_0\subseteq U$ of at least $|U|/d$ elements which has the same image in $T$. This implies  $\phi(u)^{d-1}u\ne \phi(u')^{d-1}u'$  for any $u\ne u'\in U_0$. Thus,  $|W|\ge |W_2|\ge |U_0|\ge  |U|/d$.
\end{proof}

In what follows, $d!=d(d-1)\cdots 2\cdot 1$ denotes the factorial.
\begin{corollary}\label{Schreier}
     Let $G$ be a group generated by   a   finite set $U$ as a semi-group. Let $H$ be a subgroup of index $d$ in $G$. Then there exists a finite subset $W\subset U^{\le (d!)^2-d!+1}$ such that $W$ generates $H$ as a group. Moreover, $|W|\ge \frac{|U|}{d!}$.
\end{corollary}
\begin{proof}
    Consider  the action of $G$ on the collection of left cosets of $H$ by left multiplication. This gives a  homomorphism $\rho$ from $G$ to the  permutation group $S_d$ over $d$ symbols. \revise{Thus, $\ker \rho$ is a normal subgroup of $G$ with index $[G:\ker \rho]=r\le (d!)$.}

    By Lemma \ref{SchreierNormal}, there exists a finite subset $V\subset U^{\le r^2-r+1}$ such that $V$ generates $\ker \rho$ as a group and $|V|\ge |U|/r$. As shown in Lemma \ref{SchreierNormal}, there exists a subset $T=\{a_1=1,\ldots,a_r\}\subset U^{\le r-1}$ such that $G=\bigsqcup_{i=1}^ra_i\ker \rho$. \revise{Without loss of generality, let $T\cap H=\{a_1,\ldots, a_k\}$ for $1\le k\le r$. Then we claim that $H=\bigsqcup_{i=1}^ka_i\ker \rho$. Indeed, it follows from the definition of $\rho$ and the choice of $\{a_i:1\le i\le k\}$ that $H\supset \bigsqcup_{i=1}^ka_i\ker \rho$. If there exists $k<j\le r$ such that $H\cap a_j\ker \rho\neq \varnothing$, then one gets that $a_j\in H$, which is a contradiction. Therefore, $W:=V\cup \{a_1,\ldots, a_k\}$ generates $H$ as a group. As $r\le (d!)$, one has}  $W\subset U^{\le (d!)^2-(d!)+1}$ and $|W|\ge |V|\ge |U|/{(d!)}$.
\end{proof}




\begin{lemma}\label{finiteindexPSG}
    Let $H$ be a subgroup of $G$ with index $d\geq 2$. Suppose there exist $\alpha,\beta>0$ such that for any finite set $W$ generating $H$ as a group, $|W^n|\geq (\alpha|W|)^{\beta n}$ for every $n\in \mathbb N$. Then for any finite set (possibly non-symmetric) $U$ generating $G$ as a semigroup, there exists $r=r(d)>0$ such that $$|U^n|\geq \left(\frac{\alpha}{d!2^{r/\beta}} |U|\right)^{(\beta/r)n}$$ for every $n\in \mathbb N$.
\end{lemma}
\begin{proof}
    For any finite generating set $U$ of $G$, Corollary \ref{Schreier} shows that there exists a constant $r(d)=(d!)^2-d!+1$ such that $W:=U^{\leq r}\bigcap H$ generates $H$ as a group and $|W|\ge |U|/{d!}$.

    Then for every integer $n>0$,
    $$|U^{\le n}|^r\ge |U^{\leq rn}|\geq |W^{\leq n}|\geq |W^n|\geq (\alpha|W|)^{\beta n}\geq \left(\frac{\alpha}{d!} |U|\right)^{\beta n}$$

    It follows from Remark \ref{increasing} that  $|U^n|\geq \frac{1}{n+1}|U^{\leq n}|$. Combining these two inequalities together, one has that
    $$|U^n|\ge \frac{1}{n+1}|U^{\leq n}|\ge \frac{1}{n+1}(\frac{\alpha}{d!} |U|)^{(\beta/r) n}\ge \frac{1}{2^n}\left(\frac{\alpha}{d!} |U|\right)^{(\beta/r)n}= \left(\frac{\alpha}{d!2^{r/\beta}} |U|\right)^{(\beta/r)n}$$
    for every $n\in \mathbb N$.
\end{proof}

For a symmetric version of Corollary \ref{Schreier} and Lemma \ref{finiteindexPSG}, see \cite[Prop. 2.2.26, Prop. 2.2.27]{Ker}.





\begin{lemma}\label{QuotientPSG}
   Let $\rho: G\to G'$ be an epimorphism with finite $\ker \rho$ . If there exist $\alpha,\beta>0$ such that for any finite set $U$ generating $G'$ as a semi-group (resp. group), $|U^n|\geq (\alpha|U|)^{\beta n}$ for any $n>0$. Then for any finite set $W$ generating $G$ as a semi-group (resp. group), one has $|W^n|\geq (\frac{\alpha}{|\ker \rho|} |W|)^{\beta n}$ for any $n>0$.
\end{lemma}
\begin{proof}
    \revise{As $\rho$ is an epimorphism, for any finite set $W$ generating $G$ as a semi-group (resp. group), $\rho(W)$ generates $G'$ as a semi-group (resp. group). Then one has
    $$|W^n|\ge |\rho(W)^n|\ge (\alpha |\rho(W)|)^{\beta n}\ge \left(\frac{\alpha}{|\ker \rho|}|W|\right)^{\beta n}$$ for any $n\ge 1$.}
\end{proof}


\subsection{Product set growth on one factor}\label{SSubPSGOnOneFactor}   Given an isometric group action $G\curvearrowright X$, we say that $H<G$ is a \textit{non-lineal subgroup} if  $H$   does not preserve a quasi-geodesic in $X$.
\begin{theorem}\label{PSGonOneFactor}
Let $G$ be a finitely generated group acting weakly acylindrically on a $\delta$-hyperbolic space $X$. Suppose that there exists $N>0$ so that $|E(h)\cap E(g)|\le N$ for any two independent loxodromic elements $h,g\in G$. Then there exist  constants $C, \alpha,\beta >0$ with the following property.

Let   $U\subset G$ be  a  (possibly non-symmetric) finite subset that generates a non-lineal subgroup $H$ as semigroup.  If  $\lambda (U,X)>C$, then $|U^n|\geq (\alpha|U|)^{\beta n}$ for every $n\in \mathbb N$.
\end{theorem}

The remainder of this subsection is devoted to the proof of Theorem \ref{PSGonOneFactor}. In \cite{CJY}, the product set growth for a finite symmetric $U$ was proved for non-elementary subgroups in relatively hyperbolic groups. The proof of Theorem \ref{PSGonOneFactor}  follows closely the proof  of it, in particular, the presentation of \cite[Sections 4 \& 5]{CJY}.

Briefly speaking, the main arguments in \cite{CJY} work almost word-by-word in the present setting, modulo the following two assumptions:
\begin{enumerate}
    \item
    $G$ acts weakly acylindrically on $X$, so does the subgroup $H$. This  makes through  in the current proof  the conclusion of \cite[Lemma 2.11]{CJY} (stated as Lemma \ref{LinearBILem} here), which was proved in relatively hyperbolic groups.
    \item
    There exists a uniform number $N$ so that if $g, h$ are independent, then $E(g)\cap E(h)$ has cardinality at most $N$.
\end{enumerate}
The second assumption seems to leave room for further investigation. It holds for acylindrical actions  by Lemma \ref{AcyAction}.

In what follows, we shall present the   detailed     proof of Theorem \ref{PSGonOneFactor}, and emphasize the difference along the way. We start with a brief review of general results proved in \cite[Section 4]{CJY}.

Given $c\ge 1$, let   $C=C(c,\delta)$ be the constant so that any quadrilateral with $c$-quasi-geodesic sides in a $\delta$-hyperbolic space is $C$-thin: any side is contained in the $C$-neighborhood of the other three sides.


Fix  constants $c, \theta\ge 10$.
We say that  a pair  of elements $(b,f)$ in $H$ is \textit{good} with parameters $(c,\theta)$ if there is a basepoint $o\in X$ such that the following properties (\hyperref[Property]{$\clubsuit$}) holds:
\begin{enumerate}\label{Property}
    \item[(i)]  $ b\in H$ is    a loxodromic  element with $|o-bo|>C(c,\delta)$ and $\xi:=\bigcup_{n\in \mathbb Z} b^n[o, bo]$
    is a $c$-quasi-axis of $b$.
    \item[(ii)] $ f\in H$  is an element so that $f\notin E( b)$ and $|o- fo|\le \theta |o- bo|$.
    \item[(iii)] Moreover, setting $R:=\theta |o- bo|$,  the following holds $$\text{diam}(\xi\cap N_R(f\xi))\le D\cdot ( |o- bo|+1)$$
    where   $D=D(c,\theta)$ is chosen according to Lemma \ref{LinearBILem}. 
\end{enumerate}
\begin{remark}\label{Property(iii)}
These properties  are exactly the assumptions for $(b,f)$ at the beginning of Section 4 in \cite{CJY}. As $H$ acts weakly acylindrically on $X$, property (iii) follows directly from Lemma \ref{LinearBILem}.
\end{remark}
We recall the following two results from \cite{CJY}, where    $(b,f)$ in $H$ is assumed to be  a {good} pair with parameters $(c,\theta)$.   See \cite[Figure 1]{CJY} for illustration.

\begin{lemma}\cite[Lemmas 4.1 and 4.2]{CJY}\label{Lem4.2}
    There exist constants $n_1,m_1,c_1\ge 1$ depending on $\theta,\delta$ with the following property. For any $n\ge n_1$, the element $ h:= f b^n$ is loxodromic with a $c_1$-quasi-axis $\eta$ defined as follows
    $$\eta:=\bigcup_{j\in \mathbb Z}h^j([x,y]_{\xi}[y, hx])$$ where the points $x= b^{m_1-n}o, y= b^{-m_1}o$     on the $c$-quasi-axis $\xi$ of $b$ satisfy the following inequalities
    $$
    \max\{\langle x, hx\rangle_y, \langle y, hy\rangle_{hx}\}\le C(c,\delta)
    $$
    and $d(y,hx)\ge \theta d(o,bo)$.
\end{lemma}
\begin{remark}[On the proof]\label{rmk: localquasigeod}
The above inequalities imply   that the sequence of points $$\{x,y,hx,hy, \cdots h^nx,h^ny\}$$ forms a long local quasi-geodesic.  If $\theta\ge 10$, then   $\eta$ is   a global $c_1$-quasi-geodesic (\cite[Lemma 2.4]{CJY}).
\end{remark}


\begin{lemma}\cite[Lemma 4.4]{CJY}\label{LargeInjectivity}
Under the assumption of Lemma \ref{Lem4.2}, for any $L>0$,  one may raise the value of $n$ depending on $L$ so that  the following additional properties hold for $ h= f b^n$. Denote $F=E( h)\cap E( b)$.
    \begin{enumerate}
        \item $\lambda(F,z)\le D_1$ for any $z\in \eta$, where  the constant $D_1=D_1(c_1,\delta)>0$ does not depend on $z$.
        \item for any $g\in E( h)\setminus F$, we have $\lambda(g,X)> L\cdot |o- bo|$.
        \item for any $g\in E( h)$, there exist $j\in \mathbb Z$ and $t\in F$ such that $g= h^jt$.
    \end{enumerate}
\end{lemma}

We now prepare the proof of Theorem \ref{PSGonOneFactor}.

Let $c_0$ be given by Lemma \ref{ShortHypLem} and $C_0=C(c_0,\delta)$ be the thin constant for $c_0$-quasi-geodesic quadrangles. We shall verify that $C=C_0+30\delta$ is the desired constant in Theorem \ref{PSGonOneFactor}.

\textbf{Step 1.}
Let $o\in X$ be a point satisfying $|\lambda(U,X)-\lambda(U,o)|\le \delta$. By the assumption,  $\lambda(U,X)>C=C(c_0,\delta)+30\delta$.
Let $b\in U^{\le 2}$ be  the loxodromic element  provided by Lemma \ref{ShortHypLem} so that $|o- bo|>C_0$ and $\xi:=\bigcup_{n\in \mathbb Z} b^n[o, bo]$ is a $c_0$-quasi-axis. By assumption, $H=\langle U\rangle$ does not preserve a quasi-geodesic in $X$, hence there exists $f\in U$ so that $f\notin E(b)$. Moreover, we have $|o-fo|\le 2 |o-bo|$. Indeed, Lemma \ref{ShortHypLem} gives that $$|o- bo|\ge\lambda(U, o)-11\delta\ge \frac{1}{2}\lambda(U, o).$$ As $ f\in U$, we obtain $|o- f o|\le \lambda(U, o)\le 2|o- bo|$. With Remark \ref{Property(iii)}, we have that
\begin{lemma}\label{GoodPair(bf)}
The pair $(b,f)$ is a good pair with parameters $(c_0, 2)$.
\end{lemma}
\textbf{Step 2.} We are now  using the assumption that  there exists a uniform number $N$ so that $|E(g)\cap E(h)|\le N$  for  two independent loxodromic elements $g, h$ in $H$.

Set $L:=4(m_1+1)$. Let $h:= f b^n$ be a loxodromic element given by Lemma \ref{LargeInjectivity}.  Since $ f\notin E(b)$, one has $ h=fb^n\notin E( b)$. The weakly acylindrical action $H\curvearrowright X$ implies that the axes of $h$ and $b$ have disjoint endpoints in $\partial X$, which means by definition that $ b$ and $ h$ are independent.  Hence, the assumption implies $|F|=|E(h)\cap E(b)|\le N$.

The first difference with the proof of \cite{CJY} starts from here.
Choose a largest subset $U_0$ of $U$ such that $U_0\cap F=\varnothing$ and $sF\neq s'F$ for any $s\neq s'\in U_0$. In \cite[Lemma 5.1]{CJY}, $U_0\cap F=\varnothing$ is not assumed.  As $|F|\le N$, we obtain  $|U_0|\ge |U|/N-1$. Let us assume that $|U|\ge 2N$, so $U_0$ is a non-empty set. If $|U|< 2N$, then the conclusion of Theorem \ref{PSGonOneFactor} clearly holds by choosing $\alpha=1/(2N)$ and $\beta=1$.

\begin{lemma}\label{NewVerification}
      There exists an integer $n_2\ge n_1$ depending only on $c_0,\delta$ with the following property.  For any $n\ge n_2$ and $s\in U_0$, the pair $(f b^n, s)$ is a good pair with parameters $(c_1, c_1^{-1})$.
  \end{lemma}
  \begin{proof}
      Recall that $y=b^{-m_1}o\in \eta$. By Remark \ref{Property(iii)}, it suffices for us to verify properties (i) and (ii) in (\hyperref[Property]{$\clubsuit$}) for $(fb^n,s)$ with respect to the basepoint $y$. 

  \textbf{Verification for property (i):} By Lemma \ref{Lem4.2}, $\eta$ is a $c_1$-quasi-axis of $h=fb^n$. So we only need to estimate $|y-hy|$. As $\xi:=\bigcup_{n\in \mathbb Z} b^n[o, bo]$ is a $c_0$-quasi-axis of $b$ and $|o-bo|>C_0$, one has the following estimate
  \begin{equation}\label{c0-quasi-geodesic}
      |o-b^no|\ge c_0^{-1}(n|o-bo|-c_0)> c_0^{-1}(C_0n-c_0).
  \end{equation}
  Let $C_1=C(c_1,\delta)$ be the thin constant for $c_1$-quasi-geodesic quadrangles in a $\delta$-hyperbolic space. Recall $y=b^{-m_1}o\in \eta$ and $|o-fo|\le 2|o-bo|$ from Lemma \ref{GoodPair(bf)}. Set $L=4(m_1+1)$. Let $n_2\ge n_1$ be the least integer such that
  \begin{equation}\label{ConditionI}
  \begin{array}{rl}
      |y-hy|&=|b^{-m_1}o-fb^n b^{-m_1}o|\\
      &\ge |o-b^no|-2|o-b^{m_1}o|-|o-fo|\\
      &>\max\{C_1, c_1(L|o-bo|+1+D_1)\}
  \end{array}
  \end{equation}
  holds for any $n\ge n_2$.  This completes the verification of property (i). Moreover,  (\ref{c0-quasi-geodesic}) also shows that $n_2$ only depends on $c_0,\delta$.

  \textbf{Verification for property (ii):} Recall from  the assumption that $\lambda(U,X)>C_0+30\delta$ and that $b$ is produced from Lemma \ref{ShortHypLem}.  For any $s\in U_0$, we have
 \begin{equation}\label{EstimateFors}
 \begin{array}{rl}
 |o- bo| & \ge\lambda(U, o)-11\delta \\
    &\ge \max\{\lambda(U_0, o)-11\delta,C_0+18\delta\}\\
    &\ge \max\{|o-so|-11\delta, C_0+18\delta\}.
 \end{array}
 \end{equation}
Note that $|o- b o|=| b^{-m_1}o- b b^{-m_1}o|=|y- b y|$, it follows that
\begin{align*}
    |y- s y| & = | b^{-m_1}o-s b^{-m_1}o| &\\
    & \le 2m_1|o- b o|+|o-so|  & (\text{$\Delta$-ineq.})\\
    & \le (2m_1+1)|o- b o|+11\delta  & (\ref{EstimateFors})\\
    & < L|o- b o| & (\ref{EstimateFors})\\
    & < c_1^{-1}|y-hy|-1-D_1. & (\ref{ConditionI})\\
    & < c_1^{-1}|y-hy|.
\end{align*}

On the other hand, if some $s\in U_0$ lies in $E( h)$, then Lemma \ref{LargeInjectivity} (3) gives that $s= h^jt$ for some $j\in \mathbb Z$ and $t\in F$. If $j\neq 0$, then one has
\begin{align*}
    |y- s y| & = | y-h^jty| &\\
    & \ge | y-h^jy|-|y-ty| &(\text{$\Delta$-ineq.})\\
    & \ge c_1^{-1}(|j|\cdot |y-hy|-c_1)-D_1  & (\text{Lemmas } \ref{Lem4.2}, \ref{LargeInjectivity} (1))\\
    & \ge c_1^{-1}|y-hy|-1-D_1 &  (\text{$|j|\ge 1$})\\
    & > |y-sy| &
\end{align*}
which is a contradiction. Thus, $j=0$ and then $s\in F$. However, this is  impossible, as $U_0\cap F=\varnothing$ is chosen. Hence,  the verification of property (ii) is complete.
  \end{proof}

\textbf{Step 3.}  The second   difference with \cite{CJY} lies in defining the  base for a free semi-group, instead of a free group. We have to do so, for $U$ might not be symmetric.

Set $ h= f b^{n_2}$, where $n_2$ satisfies Lemma \ref{NewVerification}. We now consider the good pairs $(h,s)$ in Lemma \ref{NewVerification} with parameters $(c_1, c_1^{-1})$, for  $s\in U_0$. Let $m_3=m_3(c_1,\delta), n_3=n_3(c_1,\delta)$ be given by Lemma \ref{Lem4.2} for each pair $(h,s)$.

\begin{lemma}
The set
\begin{equation}\label{FreeBase}
    T=\{s h^{n_3}: s\in U_0\}.
\end{equation}
generates a free semigroup of rank $|T|$ in $H$.
\end{lemma}
\begin{proof}
The point is to truncate the path labeled by  a word $W$ over $T$ (not over $T\cup T^{-1}$) in the following form
$$
W=s_1 h^{n_3}\cdot s_2 h^{n_3}\cdots s_n h^{n_3}$$
in order to produce a quasi-geodesic with the same endpoints. The proof proceeds     the same line as in proving \cite[Lemma 5.2]{CJY}. We sketch the proof and emphasize the difference.

Indeed, let us mark points $x= h^{m_3}o, y= h^{n_3-m_3}o$     on the $c_1$-quasi-axis $\eta$ of $h$.  Set $z_1=s_1x, w_1=s_1y$ and inductively,  $z_{i+1}=s_1h^{n_3}\cdots s_ih^{n_3}s_{i+1}x$ and $w_{i+1}=s_1h^{n_3}\cdots s_ih^{n_3}s_{i+1}y$ for each $i\ge 1$. See Figure \ref{fig:quasigeodesic} for   illustrating the points $z_1, w_1$ and $z_2, w_2$.  These points satisfy the following inequality
$$
\max\{\langle z_i, z_{i+1}\rangle_{w_i}, \langle w_i, w_{i+1}\rangle_{z_{i+1}}\}\le C(c_1,\delta)
$$ and $d(z_i,w_i)\ge 10 \max\{d(o,so):s\in U_0\}$.
As in Remark \ref{rmk: localquasigeod}, the sequence of points $(z_i,w_i)$ forms a chain of local quasi-geodesics, so produces a global quasi-geodesic. This proves that $T$ generates a free semigroup of rank $|T|$. We refer to \cite[Lemma 5.2]{CJY} for the remaining details.
\begin{figure}
    \centering

\tikzset{every picture/.style={line width=0.75pt}} 

\begin{tikzpicture}[x=0.75pt,y=0.75pt,yscale=-1,xscale=1]

\draw    (52.5,200) .. controls (83.5,200) and (98.5,188) .. (107.5,154) ;
\draw [shift={(107.5,154)}, rotate = 284.83] [color={rgb, 255:red, 0; green, 0; blue, 0 }  ][fill={rgb, 255:red, 0; green, 0; blue, 0 }  ][line width=0.75]      (0, 0) circle [x radius= 3.35, y radius= 3.35]   ;
\draw [shift={(52.5,200)}, rotate = 0] [color={rgb, 255:red, 0; green, 0; blue, 0 }  ][fill={rgb, 255:red, 0; green, 0; blue, 0 }  ][line width=0.75]      (0, 0) circle [x radius= 3.35, y radius= 3.35]   ;
\draw    (107.5,154) .. controls (130.5,169) and (218.5,176) .. (265.5,134) ;
\draw    (264.66,134.15) .. controls (288.66,153.77) and (307.87,153.97) .. (336.35,133.34) ;
\draw [shift={(336.35,133.34)}, rotate = 324.09] [color={rgb, 255:red, 0; green, 0; blue, 0 }  ][fill={rgb, 255:red, 0; green, 0; blue, 0 }  ][line width=0.75]      (0, 0) circle [x radius= 3.35, y radius= 3.35]   ;
\draw [shift={(264.66,134.15)}, rotate = 39.26] [color={rgb, 255:red, 0; green, 0; blue, 0 }  ][fill={rgb, 255:red, 0; green, 0; blue, 0 }  ][line width=0.75]      (0, 0) circle [x radius= 3.35, y radius= 3.35]   ;
\draw    (336.35,133.34) .. controls (344.67,159.51) and (408.37,220.62) .. (471.34,217.85) ;
\draw [shift={(471.34,217.85)}, rotate = 357.48] [color={rgb, 255:red, 0; green, 0; blue, 0 }  ][fill={rgb, 255:red, 0; green, 0; blue, 0 }  ][line width=0.75]      (0, 0) circle [x radius= 3.35, y radius= 3.35]   ;
\draw  [line width=5.25] [line join = round][line cap = round] (141.3,164.22) .. controls (141.3,164.22) and (141.3,164.22) .. (141.3,164.22) ;
\draw  [line width=5.25] [line join = round][line cap = round] (231.3,154.22) .. controls (231.3,154.22) and (231.3,154.22) .. (231.3,154.22) ;
\draw  [line width=5.25] [line join = round][line cap = round] (359.3,164.22) .. controls (359.3,164.22) and (359.3,164.22) .. (359.3,164.22) ;
\draw  [line width=5.25] [line join = round][line cap = round] (429.3,211.22) .. controls (429.3,211.22) and (429.3,211.22) .. (429.3,211.22) ;
\draw  [line width=0.75]  (142.5,161) .. controls (143.69,156.49) and (142.03,153.63) .. (137.52,152.44) -- (136.94,152.29) .. controls (130.49,150.58) and (127.87,147.47) .. (129.06,142.96) .. controls (127.87,147.47) and (124.05,148.88) .. (117.6,147.17)(120.5,147.93) -- (117.06,147.02) .. controls (112.55,145.83) and (109.69,147.49) .. (108.5,152) ;
\draw  [line width=0.75]  (262.5,133) .. controls (260.22,128.93) and (257.04,128.03) .. (252.97,130.3) -- (250.73,131.56) .. controls (244.91,134.81) and (240.86,134.4) .. (238.58,130.33) .. controls (240.86,134.4) and (239.09,138.07) .. (233.27,141.32)(235.89,139.85) -- (231.2,142.48) .. controls (227.13,144.75) and (226.23,147.93) .. (228.5,152) ;

\draw (35,197.4) node [anchor=north west][inner sep=0.75pt]    {$o$};
\draw (92,130.4) node [anchor=north west][inner sep=0.75pt]    {$s_{1} o$};
\draw (256,103.4) node [anchor=north west][inner sep=0.75pt]    {$s_{1} h^{n_{3}} o$};
\draw (108,166.4) node [anchor=north west][inner sep=0.75pt]    {$z_{1} =s_{1} x$};
\draw (208,160.4) node [anchor=north west][inner sep=0.75pt]    {$w_{1} =s_{1} y$};
\draw (319,106.4) node [anchor=north west][inner sep=0.75pt]    {$s_{1} h^{n_{3}} s_{2} o$};
\draw (453,187.4) node [anchor=north west][inner sep=0.75pt]    {$s_{1} h^{n_{3}} s_{2} h^{n_{3}} o$};
\draw (363,144.4) node [anchor=north west][inner sep=0.75pt]    {$z_{2} =s_{1} h^{n_{3}} s_2 x$};
\draw (342,208.4) node [anchor=north west][inner sep=0.75pt]    {$w_{2} =s_{1} h^{n_{3}} s_2 y$};
\draw (125,128.4) node [anchor=north west][inner sep=0.75pt]    {$h^{m_{3}}$};
\draw (221,117.4) node [anchor=north west][inner sep=0.75pt]    {$h^{m_{3}}$};

\end{tikzpicture}
    \caption{Part of the path labeled by $W$.}
    \label{fig:quasigeodesic}
\end{figure}
\end{proof}

Now, it becomes routine to complete the proof of Theorem \ref{PSGonOneFactor}. Setting $$\kappa=n_3(2n_2+1)+1, \alpha=\frac{1}{2^{\kappa+1}N},\beta=\frac{1}{\kappa},$$ we have $T\subset  U^{\leq \kappa}$.

Note that $|T|=|U_0|\geq \frac{|U|}{N}-1\ge \frac{|U|}{2N}$. The last inequality follows from the assumption that $|U|\ge 2N$. Therefore, for any $n\in \mathbb N$, one gets
$$|U^{\leq n}|^{\kappa}\ge |U^{\leq \kappa n}|\geq |T^{\leq n}|\geq |T|^n\geq \left(\frac{|U|}{N}-1\right)^n>\left(\frac{|U|}{2N}\right)^n$$
which yields
$$|U^n|\ge \frac{1}{n+1}|U^{\leq n}|\ge \frac{1}{2^n} \left(\frac{|U|}{2N}\right)^{n/\kappa}= \left(\frac{1}{2^{\kappa+1}N} |U|\right)^{n/\kappa}=(\alpha|U|)^{\beta n}.$$
This completes the proof of Theorem \ref{PSGonOneFactor}.

\subsection{Proof of Theorem \ref{PSG}}\label{SSubProofPSG}

\revise{We first  define some constants. Let $D,\delta$ be the   constants as in Theorem \ref{KeyThm}, $c$ be the constant given by Lemma \ref{LargeCommutator}, $c_0$ be the constant given by Lemma \ref{ShortHypLem}, $C=C(c_0,\delta)$ be the thin constant for $c_0$-quasi-geodesic quadrangle in a $\delta$-hyperbolic space. }

Let $U\subset G$ be any finite subset generating a group $\langle U\rangle$ as  a semi-group. If $\langle U\rangle$ is virtually abelian, then the item (1) holds; otherwise, $\langle U\rangle$ is not virtually abelian.

In the reminder of the proof,   assuming that the item (2) does not hold, we shall  show that the item (3) holds. Note that the negativity of item (2) implies that there exists a constant $N>0$ such that, for any $1\le i\le l$ and any two independent loxodromic elements $b,h$ (if there exist) on $X_i$, $|E(b)\cap E(h)\cap \langle U\rangle|\le N$.

\revise{We will subdivide our discussion in the following two cases.

\begin{itemize}
    \item \textbf{Case I: there is at least one lineal action of $\langle U\rangle$ on some factor.}

In this case, we assume that $U$ is symmetric. Let $H$ be a subgroup of $\langle U\rangle$ with index $\le 2^{l-1}$ such that the lineal actions of $H$ on each factor is orientable. According to Corollary \ref{Schreier}, there exists $r=r(l)$ such that $W:=U^{\le r}\bigcap H$ generates $H$ as a group. It follows from the symmetry of $U$ that $W$ is also symmetric.

As a result of Proposition \ref{refine}, there exists a constant $N_1=N_1(C,D,\delta,l)>0$ and $1\le i\le l$ such that $\lambda(\mathcal C(W^{\le N_1}, W^{\le N_1}), X_i)>C+30\delta$ and $H\curvearrowright X_i$ is of general type (note that $H \curvearrowright X_i$ can not be focal by Lemma \ref{alternative}). We caution the readers that the symmetry of $W$ is used to guarantee that $\mathcal C(W^{\le N_1}, W^{\le N_1})\subset W^{\le 4N_1}$.

\item \textbf{Case II: there are no lineal actions of $\langle U\rangle$ on factors.}

In this case, we set $H:=\langle U\rangle$ and $W:=U$.

As a result of Theorem \ref{KeyThm}, there exists a constant $N_2=N_2(C,D,\delta,l)>0$ such that the following holds: for any finite subset $S$ of $G$ with $|S|\ge N_2$, there exists $1\le i\le l$ such that $\lambda(S, X_i)>C+30\delta$. Since $H$ is infinite and $W$ generates $H$ as  a semi-group, $|W^{\le N_2}|=|\bigcup_{1\le k\le N_2}W^k|\ge N_2$. Hence, $\lambda(W^{\le N_2},X_i)>C+30\delta$. Then Lemma \ref{ShortHypLem} shows that the action $H\curvearrowright X_i$ has at least one loxodromic element. By assumption, the action is not lineal. Combining with Lemma \ref{alternative}, one has that $H \curvearrowright X_i$ is of general type.
\end{itemize}

In summary, for each case, we get a triple $(H,W,X_i)$ with the following properties:
\begin{itemize}
    \item[(i)] $H$ is a subgroup of $\langle U\rangle$ with index $\le 2^{l-1}$;
    \item[(ii)] $W\subset U^{\le r}$ generates $H$ as a semi-group where $r$ only depends on $l$;
    \item[(iii)] there exists a factor $X_i$ so that $\lambda(W^{\le N_0},X_i)>C+30\delta$ where $N_0:=\max\{4N_1,N_2\}$; and
    \item[(iv)] $H\curvearrowright X_i$ is of general type.
\end{itemize}

Let $\phi_i: H\to \text{Isom}(X_i)$ be the  homomorphism obtained from $\phi: H\to \prod_{i=1}^l\text{Isom}(X_i)$ post-composed with the natural $i$-th projection. As the action $H\curvearrowright X_i$ is of general type, there are at least two independent loxodromic elements $b_1,b_2$ on $X_i$. Then $\ker \phi_i\le E(b_1)\cap E(b_2)\cap H$, and thus  $|\ker \phi_i|\le N$ by the above assumption. }


Denote   $\overline H:=\phi_i(H)\leq \text{Isom}(X_i)$ and $\overline W:= \phi_i(W)$. 

To conclude the proof of Theorem \ref{PSG}, by Lemmas \ref{finiteindexPSG} and \ref{QuotientPSG}, it suffices to show that there exist constants $\alpha,\beta>0$ such that $|\overline W^n|\geq (\alpha|\overline W|)^{\beta n}$ for every $n\in \mathbb N$. This is done by applying Theorem \ref{PSGonOneFactor} to the subset $\overline W^{\le N_0}$.

\subsection{Applications}\label{SSubApps}

If  the assumption of weakly acylindrical actions on factors in Theorem \ref{PSG} are strengthened  to be acylindrical, then we can get a much simpler classification of product set growth of subgroups. \revise{This is based on the following lemma. 

\begin{lemma}\label{AcyAction}
    Let $G\curvearrowright X$ be a non-elementary acylindrical action on a hyperbolic space. Then there exists a constant $N_0>0$ such that for any two independent loxodromic elements $h,b\in G$, $|E(h)\cap E(b)|\le N_0$.
\end{lemma}
\begin{proof}
    When the action $G\curvearrowright X$ is acylindrical, \cite[Lemma 6.5]{DGO} shows that every loxodromic element $g$ is contained in a unique maximal virtually cyclic subgroup $E(g)$. This means $[E(g):\langle g\rangle]$ is finite. Then we claim that $E(h)\cap E(b)$ is a finite subgroup in $E(h)$ for any two independent loxodromic elements $h,b\in G$. Otherwise, some nontrivial power of $h$ belongs to $E(h)\cap E(b)$ which is impossible as $h$ and $b$ are independent. Moreover, \cite[Lemma 6.8]{Osi16} shows that there exists a constant $N_0>0$ such that any finite subgroup in $E(g)$ has a cardinality $\le N_0$ for any loxodromic element $g$. Thus, $|E(h)\cap E(b)|\le N_0$.
\end{proof}
}

\begin{corollary}\label{AcyPSG}
Let $G$ admit  a proper product  action on    hyperbolic spaces $X=\prod_{i=1}^lX_i$  with shadowing property on factors. Suppose that $G$ acts acylindrically on each factor. Then there exist $\alpha,\beta >0$ depending only on $\delta,l$ such that for every finite (possibly non-symmetric) subset $U\subset G$ generating a group $\langle U\rangle$ as a semi-group, one of the following  is true:
\begin{enumerate}
    \item $\langle U\rangle$ is virtually cyclic.
    \item $|U^n|\geq (\alpha|U|)^{\beta n}$ for every $n\in \mathbb N$.
\end{enumerate}
\end{corollary}
\begin{proof}
    \revise{By Remark \ref{AcyImplyWeaklyAcy} and Theorem \ref{PSG}, we need to verify:
     \begin{enumerate}
        \item[(i)] any virtually abelian subgroup must be virtually cyclic.
        \item[(ii)] there exists a constant $N_0>0$ such that for every $1\le i\le l$ and any two independent loxodromic elements (if there exist) $h,b\in G$ on $X_i$, $|E(h)\cap E(b)|\le N_0$.
        \item[(iii)] any non-virtually cyclic subgroup can not have a lineal action on a hyperbolic space.
    \end{enumerate}
   Note that the item (iii) is used to guarantee that $U$ can always be non-symmetric so long as  $\langle U\rangle$ is not virtually cyclic by the ``Moreover'' statement of Theorem \ref{PSG}.

   Let $H$ be a virtually abelian subgroup of $G$ and $H'$ be a torsion-free abelian subgroup of $H$ with finite index. Let $X_i$ be the factor space so that $H'\curvearrowright X_i$ has  a loxodromic element $h$. Since every element in $H'$ commutes with $h$, one has that $H'\le E(h)$. Recall that the action $H'\curvearrowright X_i$ is acylindrical implies that $E(h)$ is virtually cyclic \cite[Lemma 6.5]{DGO}. Hence, $H'$ is virtually cyclic and so is $H$. This proves item (i). The same arguments also show that only virtually cyclic groups can have a lineal acylindrical action on a hyperbolic space, which proves item (iii). Finally,  the item (ii) follows from Lemma \ref{AcyAction}.}
\end{proof}

Theorem \ref{AcylTreeThm} almost follows from Corollary \ref{AcyPSG} with $l=1$. However, the shadowing property needs to be replaced by  Corollary \ref{ShortHypLemOnTree}, as explained below.

\begin{proof}[Proof of Theorem \ref{AcylTreeThm}]
    In Corollary \ref{AcyPSG}, the assumption of proper product action on hyperbolic spaces with shadowing property on each factor is used only to guarantee the existence of a short loxodromic element on some factor. When $l=1$ and $X$ is a tree, Corollary \ref{ShortHypLemOnTree} gives that for any finite possibly non-symmetric $S\subset G$ without global fixed point on $X$, $S^{\le 2}$ contains a loxodromic element which has almost the same displacement of $S$. Hence, with Corollary \ref{ShortHypLemOnTree}, the proof of Corollary \ref{AcyPSG}    proves  the conclusion.
\end{proof}

In the remainder of this subsection, we describe the product set growth result for relatively hyperbolic groups, mapping class groups and 3-manifold groups.

\subsubsection*{\textbf{Relatively hyperbolic groups}}

In \cite{CJY}, Cui-Jiang-Yang obtained Theorem \ref{RHGPSGThm} under the assumption that $U$ is symmetric.  There the     symmetry assumption only occurs in \cite[Lemma 3.2]{CJY} and the final construction of free bases for  subgroups. Hence, with Lemma \ref{ShortHypLem} (a non-symmetric version of \cite[Lemma 3.2]{CJY}) and the above construction (\ref{FreeBase}) of free bases for   semi-groups, one gets the version of product set growth in Theorem \ref{RHGPSGThm}.



\subsubsection*{\textbf{Mapping class groups}}
As recalled  in Introduction,  $\mathrm{Mod}(\Sigma)$ acts on a finite product of hyperbolic spaces $\mathcal C(\mathbf Y^1)\times \cdots \times \mathcal C(\mathbf Y^l)$ so that the orbit map is a quasi-isometric embedding (see \cite{BBFa}). By \cite[Proposition 5.8]{BBFa}, there is a subgroup $G<\mathrm{Mod}(\Sigma)$ of finite index that preserves each factor $\mathcal C(\mathbf Y^i)$.

\begin{corollary}\label{MCGPSG}
    Let $\Sigma$ be an orientable finite-type surface with possibly cusps but without boundary. Then there exist $\alpha, \beta>0$ depending only on $\Sigma$ such that for every finite symmetric $U\subset \mathrm{Mod}(\Sigma)$, at least one of the following is true:
\begin{enumerate}
    \item $\langle U\rangle$ is virtually abelian.
    \item There exist a factor $\mathcal C(\mathbf Y^i)$ and two independent loxodromic elements $b,h\in \langle U\rangle$ on $\mathcal C(\mathbf Y^i)$ such that $|E(h)\cap E(b)\cap \langle U\rangle|$ is infinite.
    \item $|U^n|\geq (\alpha|U|)^{\beta n}$ for every $n\in \mathbb N$.
\end{enumerate}
Moreover, if there is no lineal action of $\langle U\rangle$ on factors, then $U$ can be non-symmetric and generates a group.
\end{corollary}
\begin{proof}
    It is known that $\mathrm{Mod}(\Sigma)$ has a finite-index pure subgroup $\Gamma_3(\Sigma)$ (see \cite[Theorem 3]{Iva}) and any pure subgroup is torsion-free.   Denote $N:=[\mathrm{Mod}(\Sigma):\Gamma_3(\Sigma)]$. Then any subgroup $H\le \mathrm{Mod}(\Sigma)$ satisfies $[H:H\cap \Gamma_3(\Sigma)]\le [\mathrm{Mod}(\Sigma): \Gamma_3(\Sigma)]=N$. As any finite pure subgroup is trivial, the above inequality implies that any finite subgroup in $\mathrm{Mod}(\Sigma)$ has a cardinality no more than $N$. This gives that: for any sequence of subgroups $\{H_n\}\subset \mathrm{Mod}(\Sigma)$, if there is no uniform bound on the cardinality of $H_n$, then there exists an infinite subgroup $H\in \{H_n\}$. 
    This turns the item (2) of Theorem \ref{PSG} into the new item (2) here.

   Shadowing property is implied in \cite[Proposition 3.10]{HP}. Moreover, \cite[Lemma 4]{TD} shows that the action on each factor is weakly-acylindrical. Hence, the conclusion just follows from Theorem \ref{PSG}.
\end{proof}

From the construction of each hyperbolic space $\mathcal C(\mathbf Y^i)$ \revise{in} \cite{BBFa}, there is no lineal action of $\mathrm{Mod}(\Sigma)$ on factors $\mathcal C(\mathbf Y^i)$ in (\ref{MCGAction}). Hence, Corollary \ref{MCGPSGCor} follows from Corollary \ref{MCGPSG}.


\subsubsection*{\textbf{3-manifold groups}}
Let  $M$ be  a connected, compact, oriented and irreducible 3-manifold. By standard 3-topology theory, $M$ has a canonical \textit{JSJ decomposition}, where the pieces  are either Seifert-fibered or admit finite-volume hyperbolic structures according to the Geometrization Conjecture. Such a 3-manifold $M$ is called a \textit{graph manifold} if all the pieces are Seifert fibered spaces; otherwise $M$ is called a \textit{mixed manifold}.

When $M$ is a mixed 3-manifold, then $\pi_1(M)$ is   hyperbolic relative to a collection of maximal graph manifold groups and  torus groups $\mathbb Z^2$ (see \cite{BW, Dah} for details).

When $M$ is a graph manifold possibly with boundary, the JSJ decomposition induces a (at most) 4-acylindrical splitting of $\pi_1(M)$ (see \cite{WZ}). Recall that a finitely generated group $G$ \textit{admits a k-acylindrical splitting} if $G$ is isomorphic to the fundamental group of a graph of groups such that the action of $G$ on the corresponding Bass-Serre tree is non-elementary and \textit{$k$-acylindrical}, that is, for any $g\neq 1\in G$, $\text{diam}(Fix(g))\le k$.


\begin{proof}[Proof of Corollary \ref{3MfdPSGCor}]

Let $H=\langle U\rangle $ be a subgroup of $\pi_1(M)$ generated by a finite possibly non-symmetric subset $U$.

Suppose first that $\pi_1(M)$ is a graph manifold group with associated Bass-Serre tree $T$. Then \revise{\cite{WZ} shows that} $H$ acts acylindrically on $T$. If $H$ does not stabilize a vertex of $T$, then it has product set growth by Theorem \ref{AcylTreeThm}. Otherwise, $H$ is a subgroup of a Seifert-fibered piece $N$ of the JSJ decomposition of $M$. Note that $\pi_1(N)$ is a central extension of a surface subgroup: it is virtually split only if the Euler number of the Seifert fiberation is zero. \revise{Recall that non-elementary hyperbolic groups have already been known to have uniform product set growth \cite{DS}.}  Thus, if $H$ fails to have product set growth, $H$ must virtually have infinite center: a finite index subgroup has infinite center.

Now we consider  a mixed manifold $M$ whose fundamental group is relatively hyperbolic. The elementary subgroups  (i.e. with at most two limit points in the Bowditch boundary) are either virtually $\mathbb Z$, virtually $\mathbb Z\times \mathbb Z$ or virtually a graph manifold group. If $H$ is elementary and has product set growth, then $H$ must be a subgroup of a graph manifold group. We are thus reduced to the above discussion. Hence, the conclusion follows.
\end{proof}

\section{Product set growth for groups acting on Hadamard manifolds}\label{SSubHadamard}
In this section, we present the proof for Theorem \ref{RiemPSGThm}. This is almost a corollary of Theorem \ref{PSG}, modulo some known facts, which we begin to introduce.

Recall that a metric space $X$ has \textit{bounded packing} with packing constant $P>0$ if every ball of radius 2 in $X$ can be covered by at most $P$ balls of radius 1.

\begin{lemma}\cite[Corollary 13.13]{BF18}\label{BP}
    For integers $P, J>0$, there exists a constant $k=k(P,J)$ with the following property. Suppose $X$ is a geodesic space and has bounded packing $P$. Let $S$ be a finite (possibly non-symmetric) set in $\text{Isom}(X)$ generating a discrete group $\langle S\rangle$ as a semi-group.
    If $\lambda(S^k, X)<J$, then $\langle S\rangle$ is virtually nilpotent.
\end{lemma}
\begin{proof}
    The proof   of \cite[Corollary 13.13]{BF18}  used a result of Breuillard-Green-Tao \cite[Corollary 11.2]{BGT} which   required the symmetry. However, the symmetry assumption in the latter only appears   in the first paragraph of the proof, where  it was only used to guarantee that $\langle S\rangle$ is a group. Hence, both results still hold when $S$ is non-symmetric and generates a group $\langle S\rangle$ as a semi-group.
\end{proof}

The following result was claimed in \cite{TD} with a sketched proof. As this fact is important in this paper, we present a detailed proof  communicated to us by T. Delzant.
\begin{lemma}\cite[Example 1]{TD}\label{WeaklyAcy}
    Let $M$ be a $r$-dimensional Riemannian manifold with curvature $-a^2\le K\le -1$ and $r\ge 2$. Then $\pi_1(M)\curvearrowright \widetilde M$ is weakly acylindrical.
\end{lemma}
\begin{proof}
    Since $\widetilde M$ is simply connected and $K\le -1$, $\widetilde M$ is CAT(-1) and thus $\delta$-hyperbolic for some $\delta>0$. Let $\mu$ be the Margulis constant of $M$. For any loxodromic isometry $g$, let $L_g$ be the geodesic axis of $g$ and $\lambda(g)$ be the minimal translation length of $g$.  As the distance function in $\mathbb H^r$ is strictly convex and the area of each geodesic quadrilateral is uniformly bounded, there exists a constant $D=D(\delta,\mu)>0$ with the following property. For any geodesic quadrilateral $[x,y][y,w][w,z][z,x] $in $ \mathbb H^r$ with $|x-z|,|y-w|\le 2\delta, |x-y|,|z-w|\ge D$, there exists a subpath $[x_1,y_1]\subset [x,y]$ with $|x_1-y_1|\ge |x-y|/D$ such that $\max_{p\in [x_1,y_1]}d(p,[z,w])\le \mu/5$.

    Let $\pi_{L_g}$ be the shortest projection map to $L_g$. Denote $$A=\{h\in \pi_1(M):  \mathrm{diam}(\pi_{L_g}(hL_g))\ge 3D(\lambda(g)+1) \}.$$ It suffices to show that $A\subset E(g)$.

    Fix an arbitrary element $h\in A$. By $\delta$-hyperbolicity of $\widetilde M$, there exists a subpath $[x,y]\subset L_g$ such that $|x-y|=3D(\lambda(g)+1)$ and $d(x,hL_g), d(y,hL_g)\le 2\delta$. As $\widetilde M$ is CAT(-1), there exists a subpath $[x_1,y_1]\subset [x,y]$ with $|x_1-y_1|\ge 3\lambda(g)$ such that $\max_{z\in [x_1,y_1]}d(z,hL_g)\le \mu/5$.

    As $|x_1-y_1|\ge 3\lambda(g)$, there exists $z\in [x_1,y_1]$ such that $[z,g^2z]=[z,gz]\cup [gz,g^2z]\subset [x_1,y_1]$. Denote $w=gz$. Then there exist $z',w'\in hL_g$ such that $|z-z'|,|w-w'|\le \mu/5$. By triangle inequality, $\big||z'-w'|-\lambda(g)\big|=\big||z'-w'|-|z-w|\big|\le |z-z'|+|w-w'|\le 2\mu/5$.

    Hence,
    \begin{equation}\label{Margulis}
        \begin{aligned}
        |z-g^{-1}hgh^{-1}z|& = |gz-hgh^{-1}z|\le |gz-w'|+|w'-hgh^{-1}z'|+|z'-z|\\ &\le |w'-hgh^{-1}z'|+2\mu/5=\big||z'-w'|-|z'-hgh^{-1}z'|\big|+2\mu/5\\ & =\big||z'-w'|-\lambda(g)\big|+2\mu/5\le 4\mu/5.
    \end{aligned}
    \end{equation}

    Denote $\pi: \widetilde M \to M$ as the covering map and $\tilde x$ as   \revise{one point so that $\pi(\tilde x)=x$}. Recall that the thin part of $M$ is defined as $M^{< \mu}:=\{x\in M: Inj_x(M)<\mu/2\}$. By thick-thin decomposition of $M$, $M^{<\mu}$ consists of a disjoint union of connected components $C_i$, each of which is homeomorphic to $\widetilde C_i/E_i$, where $\widetilde C_i$ is a connected component of $\{\tilde x\in \widetilde M: \inf_{f\neq 1\in \pi_1(M)}|\tilde x-f\tilde x|<\mu\}$ and $E_i$ is the maximal elementary subgroup which fixes $\widetilde C_i$. Actually, one has that for any element $f\in \pi_1(M)$, $f\widetilde C_i\cap \widetilde C_i\neq \emptyset$ if and only if $f\in E_i$.

    Note that (\ref{Margulis}) gives that $|z-g^{-1}hgh^{-1}z|<\mu$. If $g^{-1}hgh^{-1}=1$, then $hgh^{-1}$ fix the axis $L_g$. Hence, the axis of $hgh^{-1}$, i.e. $hL_g$, coincides with $L_g$, which implies $h\in E(g)$. If $g^{-1}hgh^{-1}\neq 1$, then (\ref{Margulis}) shows that $\pi(z)\in M^{<\mu}$. Moreover, (\ref{Margulis}) holds for any point in $[z,gz]$ and thus $\gamma=\pi([z,gz])\subset M^{<\mu}$. Denote by $C_i$ as the component containing $\gamma$, then $E_i=E(g)=\langle g_0\rangle$ where $g_0$ is a primitive hyperbolic isometry with $L_{g_0}=L_g$ and $\gamma$ is the core curve of $C_i$. Since $\pi(z)\in M^{<\mu}$, there exists $f\neq 1\in \pi_1(M)$ such that $|z-fz|=\inf_{f'\neq 1\in\pi_1(M)}|z-f'z|<\mu$. Hence, $f\widetilde C_i\cap \widetilde C_i\neq \emptyset$, which implies that $f\in E(g)=\langle g_0\rangle$. Let $f=g_0^m$ for some $m\neq 0\in \mathbb Z$. Therefore, $\lambda(g_0)=|z-g_0z|\le |z-fz|<\mu$.

    By Margulis Lemma, the set $B=\{g^{-1}hgh^{-1}: h\in A\}\cup \{g_0\}$ generates a nilpotent subgroup $K$. By \cite[Lemma 2.4]{BCG}, $K\le E(g)$. Hence, $g^{-1}hgh^{-1}\in E(g)$ and thus $h\in E(g)$.

    In conclusion, we have proved that the action $\pi_1(M)\curvearrowright \widetilde M$ is $3D$-weakly acylindrical.
\end{proof}

With Theorem \ref{PSG} and Lemmas \ref{BP} and \ref{WeaklyAcy}, we can prove Theorem \ref{RiemPSGThm} re-called below.

\begin{theorem}\label{RiemPSGThm2}
    Let $M$ be a $r$-dimensional Riemannian manifold with sectional curvature $-a^2\le K\le -1$ and $r\ge 2$. Then there exist constants $\alpha,\beta>0$ depending only on $r,a$ such that for any finite (possibly non-symmetric) $U\subset \pi_1(M)$ generating a group $\langle U\rangle$ as a semi-group, one of the following is true:
    \begin{enumerate}
    \item $\langle U\rangle$ is virtually nilpotent.
    \item $|U^n|\geq (\alpha|U|)^{\beta n}$ for every $n\in \mathbb N$.
\end{enumerate}
\end{theorem}
\begin{proof}
    Since $\widetilde M$ is simply connected and $K\le -1$, $\widetilde M$ is CAT(-1) and $\delta$-hyperbolic for some $\delta>0$. Since $-a^2\le K$ and the dimension of $M$ is $r$, we have $Ric \ge -(r-1)a^2$. As a result of \cite[Lemma 13.15]{BF18}, $\widetilde M$ has bounded packing with packing constant $Q=Q(r,a)$.

    By setting $P=Q$ and $J=[30\delta]+1$, Lemma \ref{BP} implies that there exists a constant $k=k(r,a,\delta)$ such that for any finite (possibly non-symmetric) $U\in \pi_1(M)$ generating a group $\langle U\rangle$ as  a semi-group, if $\lambda(U^k, \widetilde M)<J$, then $\Gamma:=\langle U\rangle$ is virtually nilpotent.

    Next, we consider the case that $\Gamma$ is not virtually nilpotent. Then $\lambda(U^k, \widetilde M)\ge J> 30\delta$ and $\Gamma\curvearrowright \widetilde M$ is not lineal, otherwise $\Gamma$ will be virtually cyclic by \cite[Lemma 2.5]{BCG}. By Lemma \ref{ShortHypLem}, there exist $o\in \widetilde M$ and a loxodromic element $b\in U^{\le 2k}$ with $|o-bo|\ge J-10\delta$.

    Denote by $C=C(c,\delta)$ the constant so that any quadrilateral with $c$-quasi-geodesic sides in a $\delta$-hyperbolic space is $C$-thin. Choose $f\in U$ such that $f\notin E(b)$. Up to enlarging $J$ to $J+C$ and with Lemma \ref{ShortHypLem}, the pair $(f,b)$ satisfies that: (\hyperref[Property]{$\clubsuit$})
    \begin{enumerate}
    \item[(i)] the element $b$ is loxodromic so that $|o- bo|>C$ and $\eta:=\bigcup_{n\in \mathbb Z} b^n[o,bo]$ is a $c$-quasi-axis,
    \item[(ii)] the element $f$ lies outside $E( b)$ and satisfies $|o- fo|\le 2 |o- bo|$,
    \item[(iii)] as Lemma \ref{WeaklyAcy} shows the action $\pi_1(M)\curvearrowright \widetilde M$ is weakly acylindrical, the element $ b$ satisfies the conclusion of Lemma \ref{LinearBILem}.
\end{enumerate}

    Since $\pi_1(M)$ acts freely on $\widetilde M$, then $E(g)\cap E(h)$ is trivial for any two independent loxodromic elements $g,h\in \pi_1(M)$. This implies that the second item in Theorem \ref{PSG} disappears here. Hence, with the same arguments in Theorem \ref{PSG}, we complete the proof.
\end{proof}

\bibliographystyle{amsplain}
\bibliography{UEG-Revised}

\providecommand{\bysame}{\leavevmode\hbox to3em{\hrulefill}\thinspace}
\providecommand{\MR}{\relax\ifhmode\unskip\space\fi MR }
\providecommand{\MRhref}[2]{%
  \href{http://www.ams.org/mathscinet-getitem?mr=#1}{#2}
}
\providecommand{\href}[2]{#2}
\begin{thebibliography}{10}

\bibitem{ANSGP}
Carolyn Abbott, Thomas Ng, Davide Spriano, Radhika Gupta, and Harry Petyt,
  \emph{Hierarchically hyperbolic groups and uniform exponential growth}, arXiv
  preprint arXiv:1909.00439 (2019).

\bibitem{Alp}
Roger~C. Alperin, \emph{Uniform growth of polycyclic groups}, vol.~92, 2002,
  Dedicated to John Stallings on the occasion of his 65th birthday,
  pp.~105--113. \MR{1934014}

\bibitem{AL06}
G.~N. Arzhantseva and I.~G. Lysenok, \emph{A lower bound on the growth of word
  hyperbolic groups}, J. London Math. Soc. (2) \textbf{73} (2006), no.~1,
  109--125. \MR{2197373}

\bibitem{BHS19}
Jason Behrstock, Mark Hagen, and Alessandro Sisto, \emph{Hierarchically
  hyperbolic spaces {II}: Combination theorems and the distance formula},
  Pacific Journal of Mathematics \textbf{299} (2019), no.~2, 257--338.

\bibitem{BCG}
G{\'e}rard Besson, Sylvestre Gallot, and Gilles Courtois, \emph{Uniform growth
  of groups acting on {C}artan-{H}adamard spaces}, Journal of the European
  Mathematical Society \textbf{13} (2011), no.~5, 1343--1371.

\bibitem{BBFa}
Mladen Bestvina, Ken Bromberg, and Koji Fujiwara, \emph{Constructing group
  actions on quasi-trees and applications to mapping class groups},
  Publications math{\'e}matiques de l'IH{\'E}S \textbf{122} (2015), no.~1,
  1--64.

\bibitem{BBFb}
\bysame, \emph{Proper actions on finite products of quasi-trees}, Annales Henri
  Lebesgue \textbf{4} (2021), 685--709.

\bibitem{BW}
Hadi Bigdely and Daniel~T. Wise, \emph{Quasiconvexity and relatively hyperbolic
  groups that split}, Michigan Math. J. \textbf{62} (2013), no.~2, 387--406.
  \MR{3079269}

\bibitem{BF18}
Emmanuel Breuillard and Koji Fujiwara, \emph{On the joint spectral radius for
  isometries of non-positively curved spaces and uniform growth}, Annales de
  l'Institut Fourier, vol.~71, 2021, pp.~317--391.

\bibitem{BGT}
Emmanuel Breuillard, Ben Green, and Terence Tao, \emph{The structure of
  approximate groups}, Publ. Math. Inst. Hautes \'{E}tudes Sci. \textbf{116}
  (2012), 115--221. \MR{3090256}

\bibitem{BH}
Martin~R Bridson and Andr{\'e} Haefliger, \emph{Metric spaces of non-positive
  curvature}, vol. 319, Springer Science \& Business Media, 2013.

\bibitem{BuH}
Michelle Bucher and Pierre de~la Harpe, \emph{Free products with amalgamation
  and {HNN}-extensions of uniformly exponential growth}, Matematicheskie
  Zametki \textbf{67} (2000), no.~6, 811--815.

\bibitem{BM}
Marc Burger and Shahar Mozes, \emph{Lattices in product of trees}, Publications
  Math{\'e}matiques de l'IH{\'E}S \textbf{92} (2000), 151--194.

\bibitem{But19}
J.~Button, \emph{Groups acting faithfully on trees and properly on products of
  trees}, 2019.

\bibitem{But13}
J.~O. Button, \emph{Explicit {H}elfgott type growth in free products and in
  limit groups}, J. Algebra \textbf{389} (2013), 61--77. \MR{3065992}

\bibitem{But20}
J.~O. Button, \emph{Groups acting purely loxodromically on products of
  hyperbolic graphs}, 2020.

\bibitem{But22}
J.~O. Button, \emph{Generalised {B}aumslag-{S}olitar groups and hierarchically
  hyperbolic groups}, 2022.

\bibitem{Cap}
Pierre-Emmanuel Caprace, \emph{Finite and infinite quotients of discrete and
  indiscrete groups}, Groups {S}t {A}ndrews 2017 in {B}irmingham, London Math.
  Soc. Lecture Note Ser., vol. 455, Cambridge Univ. Press, Cambridge, 2019,
  pp.~16--69. \MR{3931408}

\bibitem{CS21}
Filippo Cerocchi and Andrea Sambusetti, \emph{Entropy and finiteness of groups
  with acylindrical splittings}, Groups Geom. Dyn. \textbf{15} (2021), no.~3,
  755--799. \MR{4322011}

\bibitem{CDP}
Michel Coornaert, Thomas Delzant, and Athanase Papadopoulos,
  \emph{G{\'e}om{\'e}trie et th{\'e}orie des groupes: les groupes hyperboliques
  de gromov}, vol. 1441, Springer, 2006.

\bibitem{CS22}
R\'{e}mi Coulon and Markus Steenbock, \emph{Product set growth in {B}urnside
  groups}, J. \'{E}c. polytech. Math. \textbf{9} (2022), 463--504. \MR{4387462}

\bibitem{CJY}
Yu-Miao Cui, Yue-Ping Jiang, and Wen-Yuan Yang, \emph{Lower bound on growth of
  non-elementary subgroups in relatively hyperbolic groups}, Journal of Group
  Theory \textbf{25} (2022), no.~5, 799--822.

\bibitem{Dah}
Fran\c{c}ois Dahmani, \emph{Combination of convergence groups}, Geom. Topol.
  \textbf{7} (2003), 933--963. \MR{2026551}

\bibitem{DGO}
Fran{\c{c}}ois Dahmani, Vincent Guirardel, and Denis Osin, \emph{Hyperbolically
  embedded subgroups and rotating families in groups acting on hyperbolic
  spaces}, vol. 245, American Mathematical Society, 2017.

\bibitem{TD}
Thomas Delzant, \emph{K{\"a}hler groups, {$\mathbb R$}-trees, and holomorphic
  families of {R}iemann surfaces}, Geometric and Functional Analysis
  \textbf{26} (2016), no.~1, 160--187.

\bibitem{DS}
Thomas Delzant and Markus Steenbock, \emph{Product set growth in groups and
  hyperbolic geometry}, Journal of Topology \textbf{13} (2020), no.~3,
  1183--1215.

\bibitem{EMO}
Alex Eskin, Shahar Mozes, and Hee Oh, \emph{Uniform exponential growth for
  linear groups}, International Mathematics Research Notices \textbf{2002}
  (2002), no.~31, 1675--1683.

\bibitem{Fuj23}
Koji Fujiwara, \emph{The rates of growth in an acylindrically hyperbolic
  group}, 2023.

\bibitem{FS23}
Koji Fujiwara and Zlil Sela, \emph{The rates of growth in a hyperbolic group},
  Inventiones mathematicae (2023).

\bibitem{G}
Mikhael Gromov, \emph{Hyperbolic groups}, Springer, 1987.

\bibitem{GJN}
Radhika Gupta, Kasia Jankiewicz, and Thomas Ng, \emph{Groups acting on {CAT(0)}
  cube complexes with uniform exponential growth}, Algebraic \& Geometric
  Topology \textbf{23} (2023), no.~1, 13--42.

\bibitem{Hag21}
Mark Hagen, \emph{Non-colorable hierarchically hyperbolic groups}, Internat. J.
  Algebra Comput. \textbf{33} (2023), no.~2, 337--350. \MR{4581212}

\bibitem{HP}
Mark~F Hagen and Harry Petyt, \emph{Projection complexes and quasimedian maps},
  Algebraic \& Geometric Topology \textbf{22} (2023), no.~7, 3277--3304.

\bibitem{HTY}
Suzhen Han, Hoang~Thanh Nguyen, and Wenyuan Yang, \emph{Property {(QT)} of
  3-manifold groups}, arXiv preprint arXiv:2108.03361 (2021).

\bibitem{H}
Pierre De~La Harpe, \emph{Uniform growth in groups of exponential growth},
  Geometriae Dedicata \textbf{95} (2002), no.~1, 1--17.

\bibitem{Hug}
Sam Hughes, \emph{Lattices in a product of trees, hierarchically hyperbolic
  groups and virtual torsion-freeness}, Bull. Lond. Math. Soc. \textbf{54}
  (2022), no.~4, 1413--1419. \MR{4488315}

\bibitem{Iva}
Nikolai~V Ivanov, \emph{Subgroups of {T}eichmuller modular groups}, vol. 115,
  American Mathematical Soc., 1992.

\bibitem{KS}
Aditi Kar and Michah Sageev, \emph{Uniform exponential growth for {CAT(0)}
  square complexes}, Algebraic \& Geometric Topology \textbf{19} (2019), no.~3,
  1229--1245.

\bibitem{Ker}
Alice Kerr, \emph{Product set growth in virtual subgroups of mapping class
  groups}, arXiv preprint arXiv:2103.12643 (2021).

\bibitem{Kou}
Malik Koubi, \emph{Croissance uniforme dans les groupes hyperboliques}, Annales
  de l'institut Fourier, vol.~48, 1998, pp.~1441--1453.

\bibitem{KLN}
Robert~P. Kropholler, Rylee~Alanza Lyman, and Thomas Ng, \emph{Extensions of
  hyperbolic groups have locally uniform exponential growth}, 2020.

\bibitem{LS}
Roger~C Lyndon, Paul~E Schupp, RC~Lyndon, and PE~Schupp, \emph{Combinatorial
  group theory}, vol. 188, Springer, 1977.

\bibitem{Mah}
Johanna Mangahas, \emph{Uniform uniform exponential growth of subgroups of the
  mapping class group}, Geometric and Functional Analysis \textbf{19} (2010),
  no.~5, 1468--1480.

\bibitem{M}
Jason~Fox Manning, \emph{Quasi-actions on trees and property {(QFA)}}, Journal
  of the London Mathematical Society \textbf{73} (2006), no.~1, 84--108.

\bibitem{Osi}
Denis Osin, \emph{The entropy of solvable groups}, Ergodic Theory and Dynamical
  Systems \textbf{23} (2003), no.~3, 907--918.

\bibitem{Osi16}
\bysame, \emph{Acylindrically hyperbolic groups}, Transactions of the American
  Mathematical Society \textbf{368} (2016), no.~2, 851--888.

\bibitem{Raz}
Alexander~A. Razborov, \emph{A product theorem in free groups}, Ann. of Math.
  (2) \textbf{179} (2014), no.~2, 405--429. \MR{3152939}

\bibitem{Saf11}
S.~R. Safin, \emph{Powers of subsets of free groups}, Mat. Sb. \textbf{202}
  (2011), no.~11, 97--102. \MR{2907200}

\bibitem{ShW}
Peter~B Shalen and Philip Wagreich, \emph{Growth rates, $\mathbb z_p$-homology,
  and volumes of hyperbolic 3-manifolds}, Transactions of the American
  Mathematical Society \textbf{331} (1992), no.~2, 895--917.

\bibitem{SisHHS}
Alessandro Sisto, \emph{What is a hierarchically hyperbolic space?}, Beyond
  hyperbolicity, London Math. Soc. Lecture Note Ser., vol. 454, Cambridge Univ.
  Press, Cambridge, 2019, pp.~117--148. \MR{3966608}

\bibitem{W}
John~S Wilson, \emph{On exponential growth and uniformly exponential growth for
  groups}, Inventiones mathematicae \textbf{155} (2004), no.~2, 287--303.

\bibitem{WZ}
Henry Wilton and Pavel Zalesskii, \emph{Profinite properties of graph
  manifolds}, Geom. Dedicata \textbf{147} (2010), 29--45. \MR{2660565}

\bibitem{Wis}
Daniel~T Wise, \emph{Complete square complexes}, Commentarii Mathematici
  Helvetici \textbf{82} (2007), no.~4, 683--724.

\bibitem{Xie}
Xiangdong Xie, \emph{Growth of relatively hyperbolic groups}, Proceedings of
  the American Mathematical Society \textbf{135} (2007), no.~3, 695--704.

\end{thebibliography}
\end{document}